# Categorical and combinatorial aspects of descent theory

**Ross Street**

March 2003

**Abstract**

There is a construction which lies at the heart of descent theory. The combinatorial aspects of this paper concern the description of the construction in all dimensions. The description is achieved precisely for strict n-categories and outlined for weak n-categories. The categorical aspects concern the development of descent theory in low dimensions in order to provide a template for a theory in all dimensions. The theory involves non-abelian cohomology, stacks, torsors, homotopy, and higher-dimensional categories. Many of the ideas are scattered through the literature or are folklore; a few are new.

**Section Headings**

## 1. Introduction

Descent theory, as understood here, has been generalized from a basic example involving modules over rings. Given a ring morphism $f : R \longrightarrow S$, each right R-module M determines a right S-module $M \otimes_R S$. This process is encapsulated by the "pseudofunctor" Mod from the category of rings to the category of (large) categories; to each ring R it assigns the category ModR of right R-modules and to each morphism f the functor $- \otimes_R S : \mathrm{Mod}R \longrightarrow \mathrm{Mod}S$. The reason that Mod is not quite a functor is that the composite of ring morphisms is not taken precisely to the composite of the functors, but only up to a well-determined isomorphism. Descent data come into play when we contemplate what is needed on a right S-module N in order that it should be isomorphic to $M \otimes_R S$ for some M.

The author's interest in pseudofunctors was aroused many years ago by their appearance in group cohomology as "factor systems". It seemed inevitable that one day we would need to study even higher-dimensional weakenings of composition preservation:

up to isomorphism, then up to equivalence, and so on. Then I learned from John Roberts that cohomology itself dealt with higher-dimensional categories (where there are not only morphisms, but morphisms between morphisms — called 2-cells, and so on) and higher functors between them. He suggested studying higher-dimensional categories as the coefficient objects for non-abelian cohomology (see [Rts]). I was really taken by this idea which led to my work (see [St5]) on making a precise definition of the simplicial nerve of a strict higher category.

There are various possibilities for what we might mean by higher-dimensional categories. Initially we will concentrate on the strict ones called n-categories. While these were originally defined by Charles Ehresmann, let us recall how they were defined inductively by Eilenberg and Kelly [EK] in terms of hom enriched categories.

For any symmetric monoidal category $\mathcal{V}$, there is a symmetric monoidal category $\mathcal{V}\text{-Cat}$ whose objects are categories with homs enriched in $\mathcal{V}$; that is, $\mathcal{V}$-categories. (We will need enriched categories again later on; suitable references are [Ky] and [Bo; Chapter 6].) Starting with the category Set of sets using cartesian product for the monoidal structure, we can iterate the process $\mathcal{V} \mapsto \mathcal{V}\text{-Cat}$ yielding the following sequence of definitions:

$$\text{Set}, \quad \text{Cat} := \text{Set-Cat}, \quad 2\text{-Cat} := \text{Cat-Cat}, \quad 3\text{-Cat} := (2\text{-Cat})-\text{Cat} \quad \ldots$$

all terms having cartesian product as monoidal structure. Sets are called *0-categories*, categories are called *1-categories*, and, as we have indicated, objects of $(\text{Set-Cat})\text{-Cat}$ are called *2-categories*; and so on. Each set can be regarded as a discrete category so there are inclusions

$$\text{Set} \subset \text{Cat} \subset 2\text{-Cat} \subset 3\text{-Cat} \subset \ldots .$$

The union of this chain is the category $\omega\text{-Cat}$ of (*strict*) *ω-categories* [1].

When $\mathcal{V}$ is closed, it is enriched in itself. Each $n\text{-Cat}$ is cartesian closed and hence $n\text{-Cat}$ is itself naturally an (n+1)-category.

The *n-cells* in an ω-category can be defined recursively: the 0-cells of a set are its elements; the (n + 1)-cells of A are the n-cells of some hom n-category $A(a, b)$ for a, b objects of A. It is an important fact that n-categories are models for a finite-limit theory, in fact, a 1-sorted finite-limit theory where the one sort is "n-cell". In particular, this means that we can model n-categories in any finitely complete category $\mathcal{E}$.

Cohomology involves a "space" and a coefficients object. A fairly general notion of space is a simplicial object in some category $\mathcal{E}$. For example, in combinatorial homotopy theory, simplicial sets can act as spaces. In topos theory, the topos $\mathcal{E}$ itself is a generalized space; however, to the calculate cohomology of $\mathcal{E}$, we consider hypercovers; these are particular kinds of simplicial objects of $\mathcal{E}$.

[1] Although sometimes something a bit bigger than this union is given that name; as in [St5]. The term "∞-category" is also used.

Let $\mathbf{\Delta}$ denote the usual topologists' simplicial category; that is, the category of non-empty finite ordinals and order-preserving functions. Consider a "space" R which we consider to be a functor $\mathrm{R}:\mathbf{\Delta}^{\mathrm{op}} \longrightarrow \mathcal{C}$ of $\mathcal{C}$ (that is, a simplicial object of $\mathcal{C}$) and consider a coefficients object A which is an ω-category in $\mathcal{C}$. Form the functor $\mathcal{C}(\mathrm{R},\mathrm{A})$ : $\mathbf{\Delta} \longrightarrow \omega\text{-}\mathrm{Cat}$. We wish to construct *the cohomology ω-category* $\mathcal{H}(\mathrm{R},\mathrm{A})$ *of* R *with coefficients in* A. Some people would call this the "cocycle ω-category" rather than cohomology, but the spirit of category theory has it that our interest in cells of any ω-category is only up to the appropriate equivalence, and this very equivalence is the appropriate notion of *cobounding.*

Jack Duskin pointed out to me (probably in 1981) that the construction, called (lax) descent, should be done for any cosimplicial ω-category $\mathcal{E}:\mathbf{\Delta} \longrightarrow \omega\text{-}\mathrm{Cat}$ and should yield an ω-category $\mathrm{Desc}\,\mathcal{E}$. He proceeded to draw the diagrams for this construction in low dimensions. These diagrams are reproduced in Section 1. It then became a combinatorial challenge to make the general definition precise for all dimensions.

It was immediately clear that the objects of $\mathrm{Desc}\,\mathcal{E}$ were related closely to the "orientals" that I had introduced to define the nerve of an ω-category. The $\mathrm{n}^{\mathrm{th}}$ oriental is the "free n-category on the n-simplex". It took me quite a bit longer (surprisingly in retrospect!) to realize that the higher cells of $\mathrm{Desc}\,\mathcal{E}$ were based on the products of simplexes with "globes" (an n-globe is a "free living n-cell").

This led me to abstract the properties of simplexes that allowed the construction of free n-categories thereon. The result was the combinatorial notion of parity complex which I wanted to be closed under product. Meanwhile Michael Johnson and Robert Walters [JW] were taking a new approach to the orientals, and, in his PhD thesis, Johnson abstracted the combinatorial notion of pasting scheme. When I presented my ideas about descent and parity complexes in an Australian Category Seminar, I gave a simplistic suggestion for the product of two parity complexes. The very next week, Johnson had the correct construction. I was able to prove that parity complexes were closed under product. This involved the invention of a new order, called the "solid triangle order", on the elements of a parity complex.

We shall describe all these combinatorial matters in the present paper. We shall show how they lead to a precise definition of $\mathrm{Desc}\,\mathcal{E}$.

This paper started as a revised version of my Oberwolfach notes [St8]. However, quite a lot has happened since then. Most significantly there have been announcements of many competing definitions of weak n-category: see Leinster [Lr] for a readable discussion of most of the approaches to date. The path towards comparison of the approaches is being trod.

These developments present a further combinatorial challenge: how to construct cohomology with weak n-categories as coefficients. We shall provide some indication of an approach to this involving ideas of Batanin on computads.

Apart from the combinatorics establishing definitions and constructions, there needs to be a fully fledged theory of descent. This is worked out fairly well in what we would call dimension 2. So many of the sections of this paper are concerned with that. It is intricately related with the theory of stacks (*champs* in French) begun by Giraud [Gd].

## 2. Low-dimensional descent

In broad terms descent is about the higher categorical notion of limit. When an n-category $\mathcal{B}$ is a limit of a diagram $\mathcal{E}$ of n-categories, we can determine what data we need from the diagram $\mathcal{E}$ to "descend" to a cell of $\mathcal{B}$, uniquely up to the appropriate kind of equivalence.

For example, when $n = 0$, we know what it means for a set $\mathcal{B}$ to be the equalizer of a diagram $\mathcal{E}$ consisting of two functions $\partial_0$ and $\partial_1$ with the same domain $\mathcal{E}_0$ and codomain $\mathcal{E}_1$. An element F of $\mathcal{E}_0$ descends to a unique element of $\mathcal{B}$ if and only if $\partial_0 F = \partial_1 F$.

The example can be made slightly more complicated. Suppose we have a diagram E :

$$E_1 \overset{d_0}{\underset{d_1}{\rightrightarrows}} E_0$$

and a morphism $p : E_0 \longrightarrow B$ in a category $\mathcal{C}$. For any object X of $\mathcal{C}$, we can take the set $\mathcal{B}$ to be the homset $\mathcal{C}(B, X)$ and the functions $\partial_0$ and $\partial_1$ to be $\mathcal{C}(d_0, X)$ and $\mathcal{C}(d_1, X)$. If p exhibits $\mathcal{B}$ as the equalizer of $\mathcal{E} = \mathcal{C}(E, X)$, we may say that X *sees* B *as a coequalizer of* E. If this is true for all X, we might say that B is the *codescent object* of the diagram E. Alternatively, if E is the kernel pair of p, an X for which this is true is called a *sheaf for the cover* $p : E_0 \longrightarrow B$ of B.

Now let us look at $n = 1$. The construction of general limits of categories can be broken into various steps just as the limits of all diagrams of sets can be constructed from products and equalizers. The analogue of equalizer is what is called the *descent category* $\mathrm{Desc}\mathcal{E}$ of a diagram $\mathcal{E}$ of the form

$$\mathcal{E}_0 \;\overset{\partial_0}{\underset{\partial_1}{\rightrightarrows}}\!\!\xleftarrow{\;\iota_0\;}\; \mathcal{E}_1 \;\overset{\partial_0}{\underset{\partial_2}{\rightrightarrows}}\!\!\xrightarrow{\;\partial_1\;}\; \mathcal{E}_2$$

satisfying the usual identities for a truncated cosimplicial category:

$$\partial_s \partial_r = \partial_r \partial_{s-1} \text{ for } r < s \text{ and } \iota_0 \partial_0 = \iota_0 \partial_1.$$

The objects of $\mathrm{Desc}\mathcal{E}$ are pairs $(F, f)$ where $F$ is an object of $\mathcal{E}_0$ and $f : \partial_1 F \longrightarrow \partial_0 F$ is a morphism of $\mathcal{E}_1$ satisfying the conditions that $\iota_0 f$ is the identity morphism of $F$ and that $\partial_1 f$ is the composite of $\partial_2 f$ and $\partial_0 f$ (a commutative triangle). A morphism $u : (F, f) \longrightarrow (G, g)$ in $\mathrm{Desc}\mathcal{E}$ is a morphism $u : F \longrightarrow G$ in $\mathcal{E}_0$ such that $\partial_0 u \circ f = g \circ \partial_1 u$. Composition of morphisms in $\mathrm{Desc}\mathcal{E}$ is as in $\mathcal{E}_0$.

In particular, for categories $\mathcal{A}$ and $\mathcal{X}$, the functor category $[\mathcal{A}, \mathcal{X}]$ (whose objects are functors $\mathcal{A} \longrightarrow \mathcal{X}$ and whose morphisms are natural transformations) is the descent category for the a cosimplicial category obtained as follows. The nerve $\mathrm{Ner}\mathcal{A}$ of $\mathcal{A}$ is the simplicial set which begins

$$\mathrm{ob}\mathcal{A} \;\begin{smallmatrix} \xleftarrow{d_0} \\ \xrightarrow{i_0} \\ \xleftarrow{d_1} \end{smallmatrix}\; \mathrm{ar}\mathcal{A} \;\begin{smallmatrix} \xleftarrow{d_0} \\ \xleftarrow{d_1} \\ \xleftarrow{d_2} \end{smallmatrix}\; \mathrm{cp}\mathcal{A} \;\ldots$$

where $\mathrm{ob}\mathcal{A}$, $\mathrm{ar}\mathcal{A}$ and $\mathrm{cp}\mathcal{A}$ are the sets of objects, arrows (= morphisms), and composable pairs of arrows of $\mathcal{A}$, where the left-hand functions $d_0$, $i_0$, and $d_1$ assign codomain, identity arrow, and domain to each arrow, object, and arrow, and where the right-hand $d_1$ assigns the composite to each composable pair of arrows. We can regard each set $S$ as a discrete category; then $[S, \mathcal{X}]$ is the category of S-indexed families in $\mathcal{X}$. The cosimplicial category we want is $[\mathrm{Ner}\mathcal{A}, \mathcal{X}]$. We leave it as an exercise (although one can see a generalization in the proof of Proposition 3) to verify that there is an isomorphism of categories:

$$\mathrm{Desc}[\mathrm{Ner}\mathcal{A}, \mathcal{X}] \cong [\mathcal{A}, \mathcal{X}].$$

Because this holds naturally for all categories $\mathcal{X}$, we can re-interpret this isomorphism as saying that $\mathcal{A}$ is the *codescent category* of $\mathrm{Ner}\mathcal{A}$, showing that every category is obtained by codescent from a cosimplicial set.

The reader will need to know a little about 2-categories and (especially for the $n = 2$ case below) about pasting; an appropriate reference is [KS].

Suppose we have a truncated simplicial diagram $E$ :

$$E_2 \;\begin{smallmatrix} \xrightarrow{d_0} \\ \xrightarrow{d_1} \\ \xrightarrow{d_2} \end{smallmatrix}\; E_1 \;\begin{smallmatrix} \xrightarrow{d_0} \\ \xleftarrow{i_0} \\ \xrightarrow{d_1} \end{smallmatrix}\; E_0$$

and a morphism $p : E_0 \longrightarrow B$ in a 2-category $\mathcal{C}$. For any object $X$ of $\mathcal{C}$, we can take the category $\mathcal{B}$ to be the homset $\mathcal{C}(B, X)$ and the functions $\partial_r$ and $\iota_r$ to be $\mathcal{C}(d_r, X)$ and $\mathcal{C}(i_r, X)$. If $p$ exhibits $\mathcal{B}$ as the descent category of $\mathcal{E} = \mathcal{C}(E, X)$, we may say that $X$ *sees*

B *as a codescent object of* E. If this is true for all X, we might also say that B is the *codescent object* of the diagram E. Alternatively, if E is the simplicial kernel of p (that is, $E_1 = E_0 \times_B E_0$ is the binary product of p with itself in the slice 2-category $\mathcal{C}/B$, and $E_2 = E_0 \times_B E_0 \times_B E_0$ is the ternary product, where the morphisms $d_r$ are the projections) an X for which this is true is called a *stack for the cover* $p : E_0 \longrightarrow B$ of B.

Now let us look at n = 2. We shall describe the *descent 2-category* $\mathrm{Desc}\mathcal{E}$ of a truncated cosimplicial 2-category $\mathcal{E}$:

$$\mathcal{E}_0 \begin{array}{c} \xrightarrow{\partial_0} \\ \xleftarrow{\iota_0} \\ \xrightarrow{\partial_1} \end{array} \mathcal{E}_1 \begin{array}{c} \xrightarrow{\partial_0} \\ \xleftarrow{\iota_0} \\ \xrightarrow{\partial_1} \\ \xleftarrow{\iota_1} \\ \xrightarrow{\partial_2} \end{array} \mathcal{E}_2 \begin{array}{c} \xrightarrow{\partial_0} \\ \xrightarrow{\partial_1} \\ \xrightarrow{\partial_2} \\ \xrightarrow{\partial_3} \end{array} \mathcal{E}_3 \ .$$

The objects $(F, f, \phi)$ consist of an object F of $\mathcal{E}_0$, a morphism $f : \partial_1 F \longrightarrow \partial_0 F$ of $\mathcal{E}_1$ for which has $\iota_0 f = 1_F$, and a 2-cell

$$\begin{array}{ccc} \partial_2\partial_1 F = \partial_1\partial_1 F & \xrightarrow{\partial_1 f} & \partial_1\partial_0 F = \partial_0\partial_0 F \\ & \Downarrow \phi & \\ {\scriptstyle \partial_2 f} \searrow & \partial_2\partial_0 F = \partial_0\partial_1 F & \nearrow {\scriptstyle \partial_0 f} \end{array}$$

of $\mathcal{E}_2$ which has $\iota_0 \phi = 1_f$ and $\iota_1 \phi = 1_f$ and is such that the following equation between pasting composites holds in $\mathcal{E}_3$:

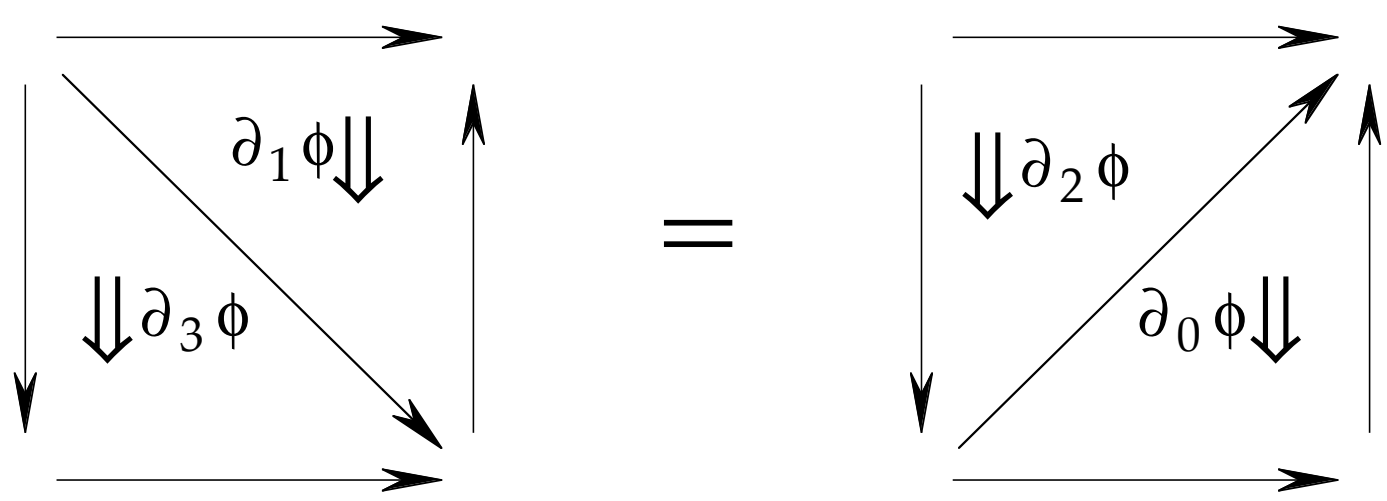

(a commutative tetrahedron). The morphisms $(u, \upsilon) : (F, f, \phi) \longrightarrow (G, g, \psi)$ consist of a morphism $u : F \longrightarrow G$ in $\mathcal{E}_0$ and a 2-cell of $\mathcal{E}_1$ which has $\iota_0\upsilon = 1_u$ and is such that the following equality (a commutative triangular cylinder) holds in $\mathcal{E}_2$.

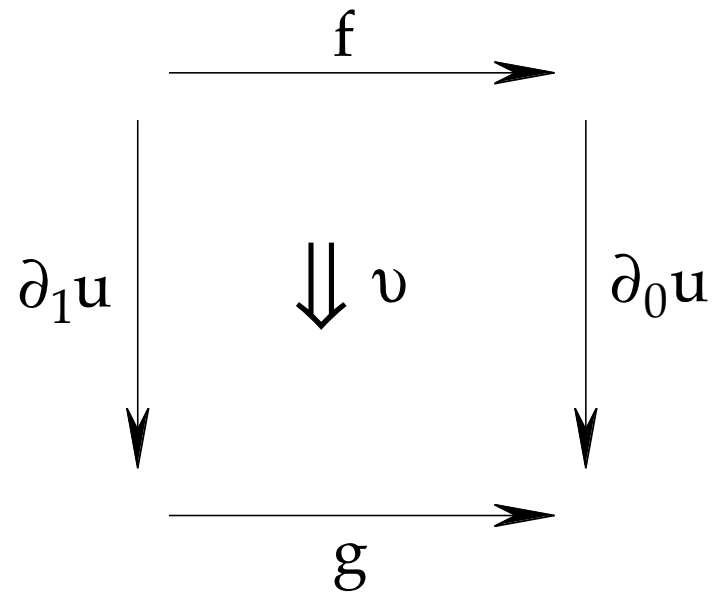

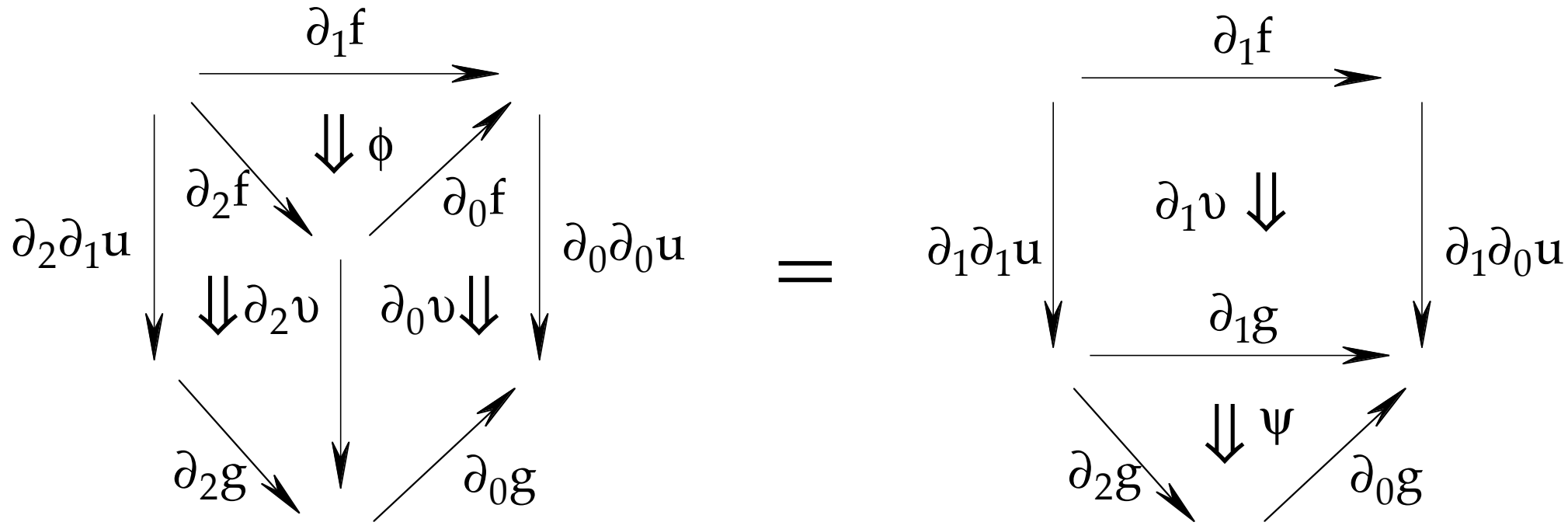

Composition of morphisms uses composition in $\mathcal{E}_0$ for the first component and vertical stacking of the 2-cells in $\mathcal{E}_1$ for the squares in the second component. The 2-cells $\alpha : (u,\upsilon) \Rightarrow (v,\nu) : (F,f,\phi) \longrightarrow (G,g,\psi)$ are just 2-cells $\alpha : u \Rightarrow v : F \longrightarrow G$ in $\mathcal{E}_0$ such that the following equality (a commutative circular cylinder) holds in $\mathcal{E}_1$.

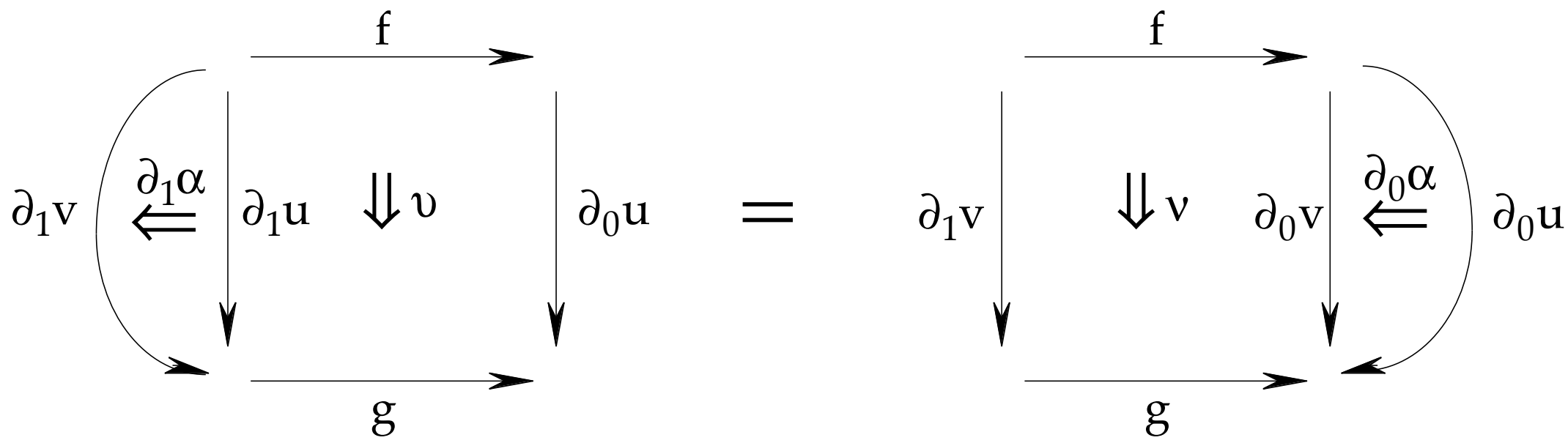

The compositions of 2-cells are those of $\mathcal{E}_0$.

Generally then, we begin with a cosimplicial $\omega$-category $\mathcal{E}$ (which is simply a functor $\mathcal{E} : \mathbf{\Delta} \longrightarrow \omega\text{-Cat}$ where $\mathbf{\Delta}$ is the *(topologists') simplicial category* whose objects are the non-empty finite ordinals) and hope to produce a *descent* $\omega$*-category* $\mathrm{Desc}\mathcal{E}$. The purpose of this paper is to make this construction precise for the case of strict $\omega$-categories, to suggest a precise construction in the case of so-called weak $\omega$-categories, and to indicate some reasons why the construction is important.

## 3. Exactness of the 2-category of categories

At the heart of modern algebra is the following exactness property of the category Set of sets. Every morphism factors, uniquely up to isomorphism, as a composite of a surjective morphism followed by an injective morphism, and the surjective morphisms are precisely those that occur as coequalizers.

When it comes to exactness properties of the 2-category Cat of categories, there are numerous possibilities. In particular, we may be interested in studying Cat as a mere bicategory in the sense of Bénabou [Bu] and work only with objects of Cat up to equivalence. However, for the moment, we wish to regard it as a (strict) 2-category and work up to isomorphism. The factorization we wish to highlight involves expressing each functor $f : A \longrightarrow B$ as a composite of functors $s : A \longrightarrow C$ and $j : C \longrightarrow B$ where $s$ is bijective on objects (b.o.) and $j$ is fully faithful (f.f.).

This defines a factorization system on the category Cat. For a given functor $f$, to produce such a factorization, define $C$ to have the objects of $A$ and the homs $C(a,a') = B(fa, fa')$; then in fact $s$ is the identity on objects and $j$ is the identity on homs. Moreover, given a commutative square of functors

$$\begin{array}{ccc} A & \xrightarrow{\;s\;} & B \\ {\scriptstyle u}\downarrow & & \downarrow{\scriptstyle v} \\ C & \xrightarrow[\;j\;]{} & D \end{array} \quad ,$$

if $s$ is b.o. and $j$ is f.f. then there exists a unique functor $w : B \longrightarrow C$ with $j\,w = v$ and $w\,s = u$; for fixed $s$ and $j$, and varying $u$ and $v$, we call this the *diagonal fill-in property*. This last property can be expressed by saying that the square

$$\begin{array}{ccc} [B,C] & \xrightarrow{\;[B,j]\;} & [B,D] \\ {\scriptstyle [s,C]}\downarrow & & \downarrow{\scriptstyle [s,D]} \\ [A,C] & \xrightarrow[\;[A,j]\;]{} & [A,D] \end{array}$$

is a pullback after applying the functor $ob : Cat \longrightarrow Set$. Actually, this last square is a pullback already in Cat. It follows that a functor $s : A \longrightarrow C$ is b.o. if and only if it has the diagonal fill-in property for all f.f. $j : C \longrightarrow B$.

On the other hand, given $s$, $j$, $u$, $v$ as above, but instead of $j\,u = v\,s$, merely an isomorphism $\sigma : j\,u \cong v\,s$, one finds that there is a unique pair $w : B \longrightarrow C$, $\tau : j\,w \cong v$ such that $w\,s = u$ and $\tau s = \sigma$. This implies that the last displayed square of functor categories is a *pseudopullback* in Cat as well as a pullback (see [JS4]).

One might feel that the b.o. functors are not the correct higher version of surjective function since a b.o. functor between discrete categories is not merely surjective but an isomorphism. This is where the bicategorical view of Cat plays its role. We call a functor $f : A \longrightarrow B$ *essentially surjective on objects* (e.s.o.) when, for all object b of B, there exists an object a of A and an isomorphism $fa \cong b$. Clearly an e.s.o. functor between discrete categories is precisely a surjective function. However, every e.s.o. functor f is equivalent to a b.o. functor; indeed, factorize $f = j\,s$ with s b.o. and j f.f., then f e.s.o. implies j is an equivalence. This means that, when regarding Cat as a bicategory, the b.o. functors are indistinguishable from the e.s.o. functors.

Now we turn to the main aspect of exactness: the higher analogue of surjective functions being coequalizers.

**Proposition 3** *A functor is bijective on objects if and only if it exhibits its codomain as the (2-categorical) codescent category of some simplicial category.*

**Proof** Suppose $p : E_0 \longrightarrow B$ exhibits B as the codescent category for a simplicial category E. This means that, for all categories X, the functor p induces an isomorphism of categories

$$[B, X] \cong \mathrm{Desc}[E, X].$$

We show that p is b.o. by showing it has the diagonal fill-in property with respect to all f.f. $j : C \longrightarrow D$. This amounts to showing that the square

$$\begin{array}{ccc} \mathrm{Desc}[E,C] & \longrightarrow & \mathrm{Desc}[E,D] \\ \downarrow & & \downarrow \\ [E_0, C] & \xrightarrow[{[E_0, j]}]{} & [E_0, D] \end{array}$$

is a pullback at the level of objects, where the vertical functors are the obvious forgetfuls. An object of $\mathrm{Desc}[E,D]$ consists of a functor $g : E_0 \longrightarrow D$ and a natural transformation $\gamma : g d_1 \longrightarrow g d_0$ satisfying conditions. So an object of the pullback consists of such data together with an object h of $[E_0, C]$ such that $g = jh$. Since j is f.f., there exists a unique natural transformation $\kappa : h d_1 \longrightarrow h d_0$ such that $j\kappa = \gamma$. Again, since j is f.f., h and $\kappa$ satisfy the descent data conditions required for h and $\kappa$ to be an object of $\mathrm{Desc}[E, C]$. So the square is indeed a pullback.

Conversely, for any functor $s : A \longrightarrow B$ there is a “higher kernel” which is a simplicial category E defined as follows. Put $E_0 = A$. Let $E_1$ be the comma category $s \downarrow s$ (in the notation of Mac Lane [ML]): the objects are triples $\left(a_1, sa_2 \xrightarrow{\beta} sa_1, a_0\right)$ where the $a_r$ are objects of A and $\beta$ is a morphism of B — the functor $d_r$ takes such a triple to $a_r$ and the functor $i_0$ takes an object a of A to $\left(a, sa \xrightarrow{1_{sa}} sa, a\right)$; the morphisms of $s \downarrow s$ are pairs of morphisms in A making the obvious square commute in B. Let $E_2$ be the category we might call $s \downarrow s \downarrow s$: the objects are quintuplets

$$\left(a_2,\ sa_2 \xrightarrow{\beta_2} sa_1, a_1, sa_2 \xrightarrow{\beta_0} sa_1, a_0\right).$$

Then there is an obvious natural transformation $\lambda : sd_1 \longrightarrow sd_0$ equipping s with the structure of an object of $\mathrm{Desc}[E, B]$. (Notice that, if s is an identity-on-objects functor from a discrete category, then E is the nerve NerB of B.)

We claim that s exhibits B as the codescent category of E if s is b.o. To see this take any category D and an object g, $\gamma : gd_1 \longrightarrow gd_0$ of $\mathrm{Desc}[E, D]$. We shall define a functor $h : B \longrightarrow D$ unique with the property that $hs = g$ and $h\lambda = \gamma$. On objects we put $hb = gs^{-1}b$. Each morphism $\beta : b_1 \longrightarrow b_0$ gives an object $e = \left(s^{-1}b_1, \beta, s^{-1}b_0\right)$ of $s \downarrow s$ and we define $h\beta = \gamma_e$. The descent conditions imply that h is indeed a functor as required. **Q.E.D.**

The final aspect of exactness of Cat that we wish to point out (making the situation much like that in a regular category in the sense of Barr [Br]) is the simple observation that *the pullback of a b.o. functor along any functor is b.o.* (After all, pullbacks in Cat are preserved by the set-of-objects functor.) There is a bicategorical analogue of this: *the pseudopullback of an e.s.o. functor along any functor is e.s.o.*

## 4. Parametrized categories

We are interested in 2-categories of categories varying over some fixed category $\mathcal{C}$. For our purposes we take a category varying over $\mathcal{C}$ to be a *pseudofunctor* $X : \mathcal{C}^{op} \longrightarrow \mathrm{Cat}$ (that is, a *homomorphism of bicategories* in the sense of [Bu]); a functor preserves composition and identities on the nose, whereas a pseudofunctor only preserves them up to coherent natural isomorphism. Between pseudofunctors there are pseudonatural transformations: these have isomorphisms in the naturality squares which satisfy the obvious coherence conditions.

We should explain a little of the folklore intuition behind such pseudofunctors. Suppose $\mathcal{C}$ is a category of sets in some universe such that $\mathcal{C}$ is actually an object of Cat. Suppose we are interested in studying categories of mathematical structures based on the sets in $C$. For example, we would be interested in the category $\mathrm{Gp}C$ of groups whose underlying sets are in $C$. Then we have no problem speaking of families of such groups parametrized by sets $\mathrm{I}$ belonging to $C$: they are merely functors $\mathrm{I} \longrightarrow \mathrm{Gp}C$. We actually have a functor $[-, \mathrm{Gp}C] : C^{\mathrm{op}} \longrightarrow \mathrm{Cat}$ taking $\mathrm{I}$ to $[\mathrm{I}, \mathrm{Gp}C]$. Notice that $[\mathrm{I}, \mathrm{Gp}C]$ is equivalent to the category $\mathrm{Gp}(C/\mathrm{I})$ of groups in the slice category $C/\mathrm{I}$. The assignment $\mathrm{I} \mapsto \mathrm{Gp}(C/\mathrm{I})$ becomes the object function of a pseudofunctor psudonaturally equivalent to $[-, \mathrm{Gp}C]$.

Suppose now that $C$ is the category of topological spaces in the universe mentioned above. One can certainly consider the category $\mathrm{Gp}C$ of topological groups (in the universe). However, in doing this, we are availing ourselves of nothing more than usual category theory. We wish to take advantage of parametrization by objects $\mathrm{I}$ of $C$. There is no obvious topology on the set of objects of $\mathrm{Gp}C$ so a functor $\mathrm{I} \longrightarrow \mathrm{Gp}C$ makes no use of topology; this time we do not have a functor $[-, \mathrm{Gp}C] : C^{\mathrm{op}} \longrightarrow \mathrm{Cat}$ available to us. A useful notion of topological group parametrized by $\mathrm{I}$ is a goup in $C/\mathrm{I}$, and we do still have a pseudofunctor $\mathrm{Gp}(C/-) : C^{\mathrm{op}} \longrightarrow \mathrm{Cat}$. In the language of parametrized category theory (in the terminology of [SS], or "indexed" category theory [PS], [Je]) over $C$, the pseudofunctor $\mathrm{Gp}(C/-)$ *is* the category of groups. We should point out here that groups give a slightly false impression of the general case since they are models of an algebraic theory — the axioms are equational. When the structures are defined using richer logic (fields or local rings, for example), it is not sufficient to take mere models in the slices $C/\mathrm{I}$.

Another good reason for looking at pseudofunctors is provided by Heller [Hr] who defines a *homotopy theory* to be a pseudofunctor $\mathrm{T} : C^{\mathrm{op}} \longrightarrow \mathrm{Cat}$ where $C$ is the category of categories in the universe we have been using above. There are some axioms on such a homotopy theory $\mathrm{T}$ including the condition that, for each morphism $\mathrm{f}$ of $C$, the functor $\mathrm{Tf}$ should have both adjoints. For example, let $\mathcal{T}$ be the category of topological spaces in the universe and define $\mathrm{TC}$ to be the homotopy category (inverting the obvious weak homotopy equivalences) of the functor category $[\mathrm{C}, \mathcal{T}]$. The adjoints of $\mathrm{Tf}$ are given by left and right homotopy Kan extensions along $\mathrm{f}$. In other words, rather than considering the mere stagnant homotopy category $\mathrm{T1}$ of $\mathcal{T}$ with its unattractive categorical properties, we consider the whole pseudofunctor $\mathrm{T}$ which, as a category parametrized by $C$, turns out to be nicely complete and cocomplete.

Let $C$ be any finitely complete category and put

$$\mathcal{F} = \mathrm{Hom}\left(C^{\mathrm{op}}, \mathrm{Cat}\right),$$

the 2-category whose objects are pseudofunctors, whose morphisms are pseudonatural transformations, and whose 2-cells are *modifications* (for example, see [KS] for precise definitions). The objects of $\mathcal{F}$ are to be thought of as large categories parametrized by $\mathcal{C}$.

A *category in* $\mathcal{C}$ can be defined to be a simplicial object $A : \mathbf{\Delta}^{op} \longrightarrow \mathcal{C}$ of $\mathcal{C}$ such that, for all objects U of $\mathcal{C}$, the simplicial set $\mathcal{C}(U, A)$ is the nerve of a category. Note that we are using the convention that $d_0 : A_1 \longrightarrow A_0$ is the *codomain* morphism and $d_1 : A_1 \longrightarrow A_0$ is the *domain* morphism for A as when we were defining the nerve of a category. A *functor between categories in* $\mathcal{C}$ is a simplicial map in $\mathcal{C}$. *Natural transformations in* $\mathcal{C}$ are defined in the obvious way yielding a 2-category $\mathrm{Cat}\mathcal{C}$ of categories in $\mathcal{C}$. Each object C of $\mathcal{C}$ gives a *discrete category in* $\mathcal{C}$; it is the constant functor $\mathbf{\Delta}^{op} \longrightarrow \mathcal{C}$ at C. In this way we regard $\mathcal{C}$ as a full subcategory of $\mathrm{Cat}\mathcal{C}$. The *opposite* $A^{op}$ of a category A in $\mathcal{C}$ is obtained by composing $A : \mathbf{\Delta}^{op} \longrightarrow \mathcal{C}$ with the functor $\mathbf{\Delta} \longrightarrow \mathbf{\Delta}$ which reverses the order on each ordinal.

Each category A in $\mathcal{C}$ gives a functor $\mathcal{C}(-, A) : \mathcal{C}^{op} \longrightarrow \mathrm{Cat}$. Any pseudofunctor pseudonaturally equivalent to such a functor $\mathcal{C}(-, A)$ is said to be an *essentially small* object of $\mathcal{F}$. This defines a Yoneda-like 2-functor $\mathrm{Cat}\mathcal{C} \longrightarrow \mathcal{F}$; since it is a fully faithful 2-functor, we identify categories in $\mathcal{C}$ with their image under it.

A Yoneda-like argument proves an equivalence of categories

$$\mathcal{F}(U, X) \simeq XU$$

which is actually pseudonatural in objects U of $\mathcal{C}$. This shows that *every pseudofunctor* X *is equivalent in the 2-category* $\mathcal{F}$ *to a 2-functor* $\mathcal{F}(-, X)$.

Given X in $\mathcal{F}$ and a category A in $\mathcal{C}$, we obtain a cosimplicial category

$$\mathbf{\Delta} \xrightarrow{\ A\ } \mathcal{C}^{op} \xrightarrow{\ \mathcal{F}(-,X)\ } \mathrm{Cat}\ .$$

Moreover, $\mathcal{F}(A, X)$ *is isomorphic to the descent category for this cosimplicial category.*

A pseudofunctor

$$\mathcal{E} : \mathbf{\Delta} \longrightarrow \mathrm{Cat}$$

might be called a *pseudocosimplicial category* : the cosimplicial identities only hold up to coherent isomorphisms. By incorporating these isomorphisms into the definition, it is possible to define a *descent category* $\mathrm{Desc}_p\mathcal{E}$ for any pseudocosimplicial category $\mathcal{E}$. Indeed, if $\mathcal{E}'$ is a cosimplicial category equivalent to $\mathcal{E}$ then there is an induced equivalence of categories

$$\mathrm{Desc}_p\mathcal{E} \simeq \mathrm{Desc}\,\mathcal{E}'.$$

This shows that, from the bicategorical point of view, the two constructions are indistinguishable when both can be made.

Returning to X in $\mathcal{F}$ and the category A in $\mathcal{C}$, we obtain a pseudocosimplicial category XA and a *generalized Yoneda-like equivalence*

$$\mathcal{F}(A,X) \simeq \mathrm{Desc}_p(XA).$$

The 2-category $\mathcal{F}$ is complete and cocomplete as a bicategory. Actually, it admits what are called pseudolimits and pseudocolimits; these can be calculated pointwise in Cat. Without going into too much detail: equalizers and pullbacks are *not* pseudolimits products, pseudopullbacks, comma categories, Eilenberg-Moore-algebra constructions, and descent categories are. For example, suppose we have morphisms

$$X \xrightarrow{f} Z \xleftarrow{g} Y$$

in $\mathcal{F}$. We can form both the *comma object* $f \downarrow g$ and the *pseudopullback* P of f and g as objects of $\mathcal{F}$; it is done componentwise:

$$(f \downarrow g)U = f_U \downarrow g_U$$

and PU is the full subcategory consisting of the objects

$$\left(x \in XU,\ f_U x \xrightarrow{\zeta} g_U y,\ y \in YU\right)$$

with $\zeta$ invertible. Because in these definitions we are not asking any objects to actually be equal, both $f \downarrow g$ and P can be defined on morphisms of $\mathcal{C}$ making them objects of $\mathcal{F}$. There are pseudonatural projections $p : f \downarrow g \longrightarrow X$ and $q : f \downarrow g \longrightarrow Y$ with component at U taking $(x,\zeta,y)$ to x and y, and a modification $\lambda : f\,p \longrightarrow g\,q$ with component at U having component $\zeta$ at $(x,\zeta,y)$.

Suppose X is an object of $\mathcal{F}$ and we have an object x of XU and an object y of XV which we identify with morphisms $x : U \longrightarrow X$ and $y : V \longrightarrow X$ in $\mathcal{F}$ (with U and V in $\mathcal{C}$). The *hom of* x *and* y is the comma object $x \downarrow y$. We say the hom is *small* when $x \downarrow y$ is essentially small. In this case there is a span $U \xleftarrow{p} X(x,y) \xrightarrow{q} V$ in $\mathcal{C}$ which is equivalent in $\mathcal{F}$ to $U \xleftarrow{p} x \downarrow y \xrightarrow{q} V$. We call $x : U \longrightarrow X$ *(left) homly* when the hom of x and y is small for all $y : V \longrightarrow X$. (I have used the word "admissible" in the past but this was met with objections!) A morphism $f : Z \longrightarrow X$ in $\mathcal{F}$ is called *homly* when, for all $z : U \longrightarrow Z$ with $U$ in $\mathcal{C}$, the composite $fz : U \longrightarrow X$ is homly.

A functor $p : E \longrightarrow A$ between categories in $\mathcal{C}$ is said to be a *discrete fibration* when the commutative square

$$\begin{array}{ccc} E_1 & \xrightarrow{p_1} & A_1 \\ {\scriptstyle d_1}\downarrow & & \downarrow{\scriptstyle d_1} \\ E_0 & \xrightarrow[p_0]{} & A_0 \end{array}$$

is a pullback in $\mathcal{C}$. The composite of two discrete fibrations is a discrete fibration. A discrete fibration into a discrete category has discrete domain; every morphism of $\mathcal{C}$ is such. For any functor $f : B \longrightarrow A$ in $\mathcal{C}$, the pullback of a discrete fibration $p : E \longrightarrow A$ along $f$ is a discrete fibration $p_f : E_f \longrightarrow B$. A functor $q : E \longrightarrow B$ between categories in $\mathcal{C}$ is called a *discrete opfibration* when $q^{op} : E^{op} \longrightarrow B^{op}$ is a discrete fibration.

There is a two-sided version of discrete fibration. A span

$$A \xleftarrow{\;p\;} E \xrightarrow{\;q\;} B$$

in $\mathrm{Cat}\mathcal{C}$ is called a *discrete fibration from* $A$ *to* $B$ when, in the diagram below, where the diamonds are pullbacks, $p\,i_r : E_r \longrightarrow A$ is a discrete fibration and $q\,i_\ell : E_\ell \longrightarrow B$ is a discrete opfibration.

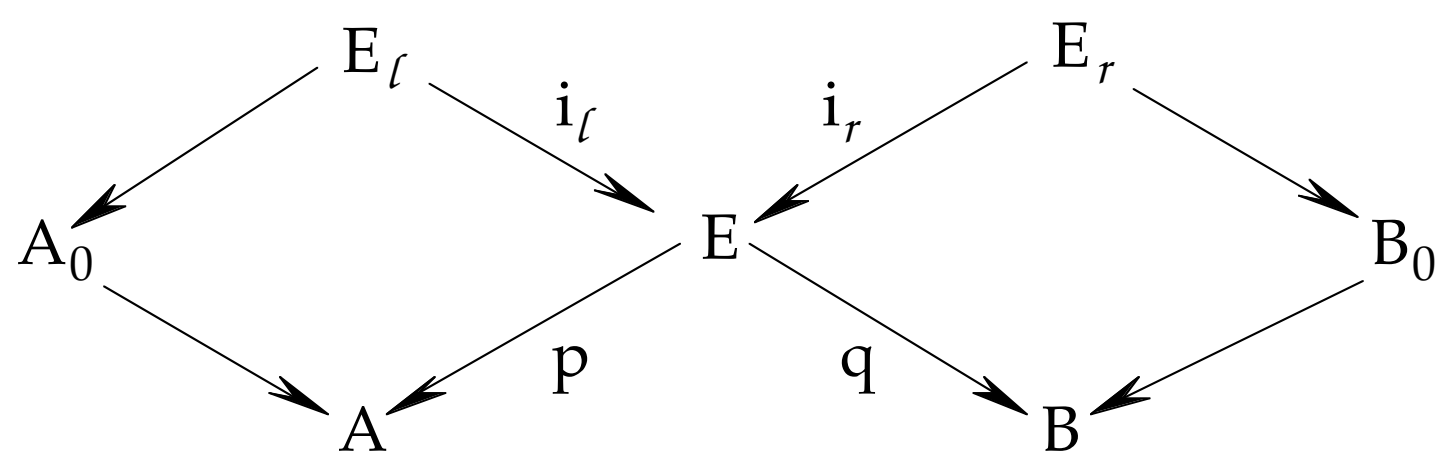

When $B$ is discrete, this reduces to the requirement that $p$ should be a discrete fibration. For all functors $u : A \longrightarrow C$ and $v : B \longrightarrow C$ in $\mathcal{C}$, the span $A \xleftarrow{\;p\;} u \downarrow v \xrightarrow{\;q\;} B$ is a discrete fibration from $A$ to $B$.

Let $A$ be a category in $\mathcal{C}$. We shall define an object $\mathcal{P}A$ of $\mathcal{F}$ called the *presheaf object of* $A$. For each object $U$ of $\mathcal{C}$, the category $(\mathcal{P}A)U$ has as objects discrete fibrations $(p, E, q)$ from $A$ to $U$ (sometimes written in abbreviated notation as $E$). We make $\mathcal{P}A$ pseudofunctorial by using pullback: that is, for $u : V \longrightarrow U$, define $(\mathcal{P}A)u$ to take $(p, E, q)$ to $(pu', E', q')$ where $E'$ is the pullback of $q$ and $u$ with projections $u'$ and $q'$.

There is a *yoneda morphism* $y_A : A \longrightarrow \mathcal{P}A$ whose component at $U$ takes $a \in AU$ to the span

$$A \xleftarrow{\ p\ } A\downarrow a \xrightarrow{\ q\ } U$$

obtained as the comma object of $1_A : A \longrightarrow A$ and $a : U \longrightarrow A$ in $\mathrm{Cat}\mathcal{C}$. The fact that the yoneda morphism is f.f. follows from the following "Yoneda lemma" in $\mathrm{Cat}\mathcal{C}$.

**Lemma 4.1** *Suppose* $a : B \longrightarrow A$ *is a functor between categories in* $\mathcal{C}$ *and* $E$ *is a discrete fibration from* $A$ *to* $B$. *Then there is an isomorphism between the category of span morphisms from* $A\downarrow a$ *to* $E$ *and the category of span morphisms from* $(a, B, 1_B)$ *to* $E$. *The isomorphism is given by composing with the right adjoint* $i : B \longrightarrow A\downarrow a$ *of* $q$.

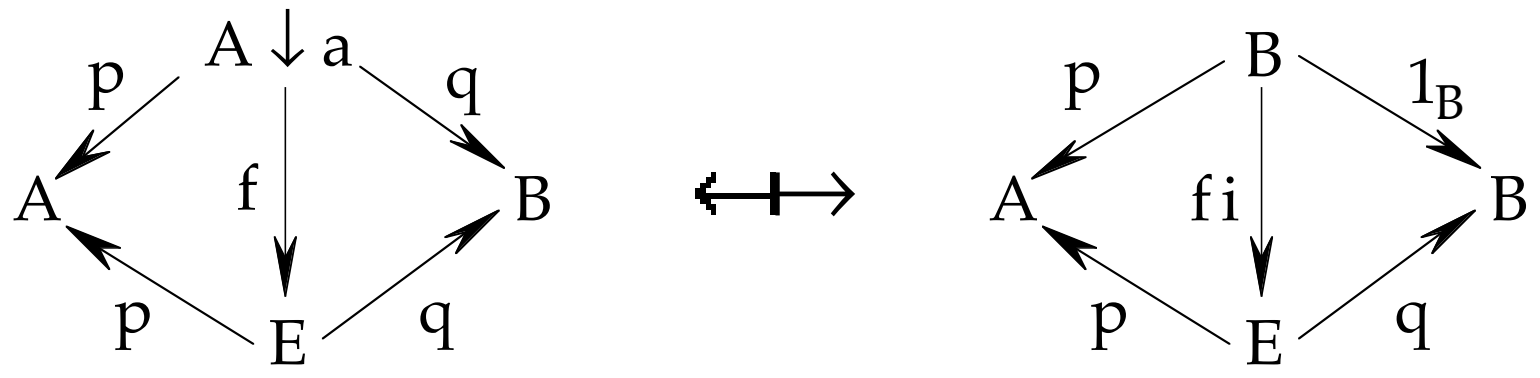

A morphism $u : W \longrightarrow U$ in $\mathcal{C}$ is said to be *powerful* (or "exponentiable") when the functor $\mathcal{C}/U \longrightarrow \mathcal{C}$, taking $K \longrightarrow U$ to the pullback $K \times_U W$, has a right adjoint. This is equivalent to asking that the functor $\mathcal{C}/U \longrightarrow \mathcal{C}/W$, taking $K \longrightarrow U$ to $K \times_U W \longrightarrow W$, have a right adjoint. It is also equivalent to the requirement that the functor $\mathcal{C}/U \longrightarrow \mathcal{C}/U$, defined by taking binary product with the object $u$ of $\mathcal{C}/U$, should have a right adjoint (so that $u$ can be used as a *power* for cartesian exponentiation in the category $\mathcal{C}/U$). Any pullback of a powerful morphism is powerful and any composite of powerful morphisms is powerful. Every morphism in a topos is powerful and the powerful morphisms in Cat were characterized by Giraud [Gd] and Conduché [Cé]; it was extended to categories in a topos by Johnstone [Je].

To see the relevance of the powerful morphisms, consider the special case of $\mathcal{P}A$ when $A$ is the terminal object $1$ of $\mathcal{C}$. Then $(\mathcal{P}1)U$ is the slice category $\mathcal{C}/U$. Morphisms $x : U \longrightarrow \mathcal{P}1$ and $y : V \longrightarrow \mathcal{P}1$ in $\mathcal{F}$ can be identified with morphisms $x : S \longrightarrow U$ and $y : S \longrightarrow V$ in $\mathcal{C}$. If $x : S \longrightarrow U$ is powerful then so too is the morphism $x \times V : S \times V \longrightarrow U \times V$; so the internal hom of the objects

$$x \times V : S \times V \longrightarrow U \times V \quad \text{and} \quad U \times y : U \times T \longrightarrow U \times V$$

of $\mathcal{C}/U \times V$ exists. This provides the span $U \xleftarrow{\ p\ } (\mathcal{P}1)(x, y) \xrightarrow{\ q\ } V$ as in the definition of small homs. It follows that the powerful morphisms $x : S \longrightarrow U$ are the homly objects of $(\mathcal{P}1)U$. With a little more work one can show that:

**Proposition 4.2** *If* A *is a category in* $\mathcal{C}$ *and* $(p, E, q)$ *is an object of* $(\mathcal{P}A)U$ *for which* $q_0 : E_0 \longrightarrow U$ *is a powerful morphism of* $\mathcal{C}$ *then* $(p, E, q)$ *is homly.*

**Corollary 4.3** *The yoneda morphism of* $A$ *is homly if* $d_0 : A_1 \longrightarrow A_0$ *is powerful.*

Recall that, in any bicategory $\mathcal{K}$, a diagram

$$\begin{array}{ccc} A & \xrightarrow{\quad j \quad} & B \\ & {\scriptstyle h}\searrow \overset{\kappa}{\Rightarrow} \swarrow{\scriptstyle k} & \\ & X & \end{array}$$

is said to exhibit k as a *left extension* of h along j when, for all morphisms $g : B \longrightarrow X$, the function

$$\mathcal{K}(B,X)(k,g) \longrightarrow \mathcal{K}(A,X)(h,g\,j)\,, \quad \sigma \mapsto \sigma j \circ \kappa$$

is a bijection. Assume $\mathcal{K}$ admits all comma objects. The diagram is said to exhibit k as a pointwise left extension of h along j when, for all $s : C \longrightarrow B$, the diagram

$$\begin{array}{ccc} j \downarrow s & \xrightarrow{\quad q \quad} & C \\ {\scriptstyle p}\downarrow & \overset{\lambda}{\Rightarrow} & \downarrow{\scriptstyle s} \\ A & \xrightarrow{\quad j \quad} & B \\ & {\scriptstyle h}\searrow \overset{\kappa}{\Rightarrow} \swarrow{\scriptstyle k} & \\ & X & \end{array}$$

exhibits $k\,s$ as a left extension of $h\,p$ along $q$. It can be shown that pointwise left extensions are indeed left extensions. For $\mathcal{K}$ = Cat, a left extension of h along a functor into 1 gives the colimit of h; and Lawvere's colimit formula for the left Kan extension is obtained from the pointwise condition with $C = 1$. For $\mathcal{K} = \mathcal{F}$, to test whether a left extension is pointwise, it suffices to take C in $\mathcal{C}$.

Suppose A is a category in $\mathcal{C}$ and $E = (p, E, q)$ is an object of $(\mathcal{P}A)U$. Let $f : A \longrightarrow X$ be a morphism in $\mathcal{F}$. The *colimit of* f *weighted by* E is the pointwise left extension $\mathrm{col}(E, f) : U \longrightarrow X$ of $f\,p$ along $q$ (see the diagram below); sometimes we identify $\mathrm{col}(E, f)$ with the corresponding object of XU.

$$\begin{array}{ccc} E & \xrightarrow{\;q\;} & U \\ {\scriptstyle p}\downarrow & \overset{\kappa}{\Rightarrow} & \downarrow{\scriptstyle \mathrm{col}(E,f)} \\ A & \xrightarrow[\;f\;]{} & X \end{array}$$

We call an object X of $\mathcal{F}$ *(small) cocomplete* when col(E,f) exists for all categories A in $\mathcal{C}$, all U in $\mathcal{C}$, all E in $(\mathcal{P}A)U$, and all $f : A \longrightarrow X$.

**Proposition 4.4** *For all categories* B *in* $\mathcal{C}$, *the presheaf object* $\mathcal{P}B$ *is small cocomplete.*

**Proof** The category of discrete fibrations from B to A is equivalent to $\mathcal{F}(A,\mathcal{P}B)$; given a discrete fibration $(p,F,q)$ from B to A, there is a corresponding $f : A \longrightarrow \mathcal{P}B$ whose component at U takes each $a : U \longrightarrow A$ to the discrete fibration from B to U defined by pulling back q and a. To obtain the colimit of f weighted by a discrete fibration E from A to U, one merely composes F from B to A with E from A to U to obtain $E \circ F$ from B to U. Then $\mathrm{col}(E,f) = E \circ F$ in $\mathcal{P}B$. **Q.E.D.**

Suppose $u : V \longrightarrow U$ is a morphism of $\mathcal{C}$ and X is an object of $\mathcal{F}$. If the functor $Xu : XU \longrightarrow XV$ has a left adjoint $\hat{X}u : XV \longrightarrow XU$ then every morphism $y : V \longrightarrow X$ has a left extension along u ; to calculate it, we use the Yoneda-like correspondence to identify y with an object of XV — then the left extension is the morphism $U \longrightarrow X$ identified with the object $(\hat{X}u)y$ of XU. Conversely, if the left extension exists for all y then Xu has a left adjoint. Pointwiseness of the left extension is equivalent to a so-called *Beck-Chevalley condition*: for every pullback

$$\begin{array}{ccc} P & \xrightarrow{\;q\;} & W \\ {\scriptstyle p}\downarrow & & \downarrow{\scriptstyle s} \\ V & \xrightarrow[\;u\;]{} & U \end{array}$$

in $\mathcal{C}$, the functor Xq has a left adjoint $\hat{X}q$ and the mate

$$\begin{array}{ccc} XP & \xrightarrow{\;\hat{X}q\;} & XW \\ {\scriptstyle Xp}\uparrow & \Rightarrow & \uparrow{\scriptstyle Xs} \\ XV & \xrightarrow[\;\hat{X}u\;]{} & XU \end{array}$$

of the canonical isomorphism $Xp \, . \, Xu \cong Xq \, . \, Xs$ is invertible.

In fact, to have $X$ such that each $Xu : XU \longrightarrow XV$ has a left adjoint and the Beck-Chevalley condition holds is equivalent to having all pointwise left extensions existing for morphisms into $X$ along morphisms in $\mathcal{C}$. It is reasonable to say in this case that $X$ has *small coproducts.* (Compare the case of Cat where all Kan extensions into $X$ exist along functors between small discrete categories if and only if $X$ admits all small coproducts.)

**Proposition 4.5** *If* $X : \mathcal{C}^{op} \longrightarrow \mathrm{Cat}$ *is a functor that has small coproducts as an object of* $\mathcal{F}$ *then, for every category* $A$ *in* $\mathcal{C}$, *the forgetful functor* $\mathrm{Desc}XA \longrightarrow XA_0$ *is monadic. The underlying functor of the monad is the composite*

$$XA_0 \xrightarrow{Xd_1} XA_1 \xrightarrow{\hat{X}d_0} XA_0 .$$

**Proof** The equation $d_0 i_0 = 1$ induces a natural transformation $Xi_0 \Rightarrow \hat{X}d_0$ which then restricts along $Xd_1$ to give a natural transformation $\eta : 1_{XA_0} \Rightarrow \hat{X}d_0 . Xd_1$ (where we use $d_1 i_0 = 1$). Using the units of the adjunctions $\hat{X}d_1 \dashv Xd_1$ and then those of the adjunctions $\hat{X}d_0 \dashv Xd_0$ and $\hat{X}d_2 \dashv Xd_2$, we obtain a natural transformation

$$1_{XA_2} \Rightarrow Xd_1 . Xd_0 . \hat{X}d_0 . Xd_1 . \hat{X}d_1 . \hat{X}d_1 .$$

Using the equations $d_0 d_1 = d_0 d_0$ and $d_1 d_1 = d_1 d_2$, we see that the codomain of this natural transformation is isomorphic to $Xd_0 . Xd_0 . \hat{X}d_0 . Xd_1 . \hat{X}d_1 . \hat{X}d_2$ and so, using mates under the adjunctions $\hat{X}d_0 \dashv Xd_0$ and $\hat{X}d_2 \dashv Xd_2$, we obtain a natural

$$\hat{X}d_0 . Xd_2 \Rightarrow Xd_0 . \hat{X}d_0 . Xd_1 . \hat{X}d_1 .$$

Since $d_0$ and $d_2$ exhibit $A_2$ as a pullback of $d_0$ and $d_1$, the Beck-Chevalley condition gives $\hat{X}d_0 . Xd_2 \cong Xd_1 . \hat{X}d_0$. So we have a natural transformation

$$Xd_1 . \hat{X}d_0 \Rightarrow Xd_0 . \hat{X}d_0 . Xd_1 . \hat{X}d_1$$

which has a mate $\mu : \hat{X}d_0 . Xd_1 . \hat{X}d_0 . Xd_1 \Rightarrow \hat{X}d_0 . Xd_1$. The functor $\hat{X}d_0 . Xd_1$ becomes a monad on the category $XA_0$ by taking $\eta$ and $\mu$ as unit and multiplication. An Eilenberg-Moore algebra for this monad is an object $x$ of $XA_0$ together with an action $(\hat{X}d_0 . Xd_1)x \longrightarrow x$ which corresponds under the adjunction $\hat{X}d_0 \dashv Xd_0$ to a morphism $\xi : (Xd_1)x \longrightarrow (Xd_0)x$; the action conditions translate to the conditions that $(x,\xi)$ should be an object of DescXA. This defines an isomorphism between the category of Eilenberg-Moore algebras and DescXA. **Q.E.D.**

## 5. Factorizations for parametrized functors

We shall make use of the exactness properties of Cat (see Section 3) carried over, in a pointwise manner to $\mathcal{F} = \mathrm{Hom}(\mathcal{C}^{op}, \mathrm{Cat})$. Consider first the factorization into b.o. and f.f.

Let $f : X \longrightarrow Y$ be any morphism of $\mathcal{F}$. We can factorize each component functor $f_U : XU \longrightarrow YU$ into a composite of a b.o. $s_U : XU \longrightarrow ZU$ and f.f. $j_U : ZU \longrightarrow YU$. We would like to make $Z$ into an object of $\mathcal{F}$; that is, a pseudofunctor. For all $u : V \longrightarrow U$ in $\mathcal{C}$, we have an isomorphism

$$\begin{array}{ccccc} XU & \xrightarrow{s_U} & ZU & \xrightarrow{j_U} & YU \\ {\scriptstyle Xu}\downarrow & & \cong \; f_u & & \downarrow{\scriptstyle Yu} \\ XV & \xrightarrow[s_V]{} & ZV & \xrightarrow[j_V]{} & YV \end{array}$$

which, by the 2-categorical diagonal fill-in property, is equal to

$$\begin{array}{ccccc} XU & \xrightarrow{s_U} & ZU & \xrightarrow{j_U} & YU \\ {\scriptstyle Xu}\downarrow & & {\scriptstyle Zu}\downarrow & \cong \; j_u & \downarrow{\scriptstyle Yu} \\ XV & \xrightarrow[s_V]{} & ZV & \xrightarrow[j_V]{} & YV \end{array}$$

for a unique functor $Zu$ and isomorphism $j_u$, where the left square commutes. Using the uniqueness of this kind of fill-in, we see that $Z$ becomes a pseudofunctor and that $s : X \longrightarrow Z$ and $j : Z \longrightarrow Y$ become pseudonatural; in fact, each $s_u$ is an identity (that is, it is *strict*).

In summary, every morphism $f$ of $\mathcal{F}$ has the form $f = j\,s$ where $s$ is pointwise b.o. and strict, and $j$ is pointwise f.f. From the bicategorical view, this is much stricter than we need. Let us recall the bicategorical notion of factorization system.

Suppose $\mathcal{M}$ is a bicategory with two distinguished classes $\mathcal{S}$ and $\mathcal{J}$ of morphisms. We call the pair $(\mathcal{S}, \mathcal{J})$ a *factorization system on the bicategory* $\mathcal{M}$ when

(i) each of $\mathcal{S}$ and $\mathcal{J}$ contains the equivalences and is closed under composition;

(ii) for each morphism $f$ of $\mathcal{M}$, there exist $s$ in $\mathcal{S}$, $j$ in $\mathcal{J}$ and an isomorphism
$$f \cong j\,s;$$

(iii) for each $s : X \longrightarrow K$ in $\mathcal{S}$ and $j : Z \longrightarrow Y$ in $\mathcal{J}$, the following square is equivalent to the pseudopullback of the left and bottom sides.

$$\begin{array}{ccc} \mathcal{M}(K,Z) & \xrightarrow{\mathcal{M}(K,j)} & \mathcal{M}(K,Y) \\ {\scriptstyle \mathcal{M}(s,Z)}\downarrow & \cong & \downarrow{\scriptstyle \mathcal{M}(s,Y)} \\ \mathcal{M}(X,Z) & \xrightarrow[\mathcal{M}(X,j)]{} & \mathcal{M}(X,Y) \end{array}$$

The factorization system is called *regular* when the pseudopullback, along any morphism, of a morphism in $\mathcal{S}$ is also in $\mathcal{S}$.

For our pointwise b.o./f.f. factorization system on $\mathcal{F}$, condition (i) already causes a problem since equivalences are not necessarily b.o. It is necessary (as implied in Section 3) to allow the more general pointwise e.s.o. morphisms in place of the pointwise b.o. Let us call morphisms of $\mathcal{F}$ *b.o., e.s.o.,* or *f.f.* when they are pointwise so.

The analysis of Section 3 and the remarks at the beginning of this section make it easy to see that *the classes of e.s.o. and f.f. morphisms form a regular factorization system on $\mathcal{F}$ as a bicategory.* Notice that the pointwiseness of the morphism classes can be expressed by the fact that, for U in $\mathcal{C}$, the 2-functor

$$\mathcal{F}(U,-) : \mathcal{F} \longrightarrow \mathrm{Cat}$$

preserves the bicategorical e.s.o./f.f. factorization (we use the Yoneda-like equivalence between this 2-functor and evaluation at U). Moreover, the f.f. morphisms in $\mathcal{F}$ can be characterized as those morphisms $j : Z \longrightarrow Y$ for which the functors

$$\mathcal{F}(X,j) : \mathcal{F}(X,Z) \longrightarrow \mathcal{F}(X,Y)$$

are fully faithful for all X in $\mathcal{F}$. The e.s.o. morphisms are not preserved by all $\mathcal{F}(X,-)$, however, they do enjoy the codescent characterization:

**Proposition 5.1** *A morphism of $\mathcal{F}$ is (pointwise) essentially surjective on objects if and only if it exhibits its codomain as the bicategorical codescent category of some pseudosimplicial object of $\mathcal{F}$.*

**Proof** By using the factorization described at the beginning of this section, we can factor each e.s.o. in $\mathcal{F}$ into a pointwise b.o. followed by an equivalence. Since we are only interested in $\mathcal{F}$ as a bicategory, we can work with this pointwise b.o. $s : X \longrightarrow Z$ (which can also be assumed strict). Now we can use pointwise the "generalized kernel" construction in the proof of Proposition 3. Notice that the construction involves comma categories and the like which are pseudolimits and so create a simplicial object of $\mathcal{F}$. Since the descent construction is a pseudolimit too, the codescent construction of our simplicial object of $\mathcal{F}$ is formed pointwise. It therefore follows from Proposition 3 that s induces an isomorphism of Z with the codescent object.

The converse also follows from Proposition 3 and the pointwise nature of the codescent construction. **Q.E.D.**

Suppose now that the category $\mathcal{C}$ has a class of distinguished morphisms called *covers*. We assume that covers form a calculus of left fractions; this means, they contain the isomorphisms, are closed under composition, and, for each object U, the opposite of the full subcategory CovU of $\mathcal{C} \downarrow U$, consisting of the covers, is filtered. A trivial example is when the covers are precisely the isomorphisms. A more interesting example is when $\mathcal{C}$ is a regular category and the covers are the strong epimorphisms (which are the same as the extremal and regular epimorphisms for $\mathcal{C}$ regular).

A morphism $f : X \longrightarrow Y$ in $\mathcal{F}$ is called *locally surjective on objects (l.s.o.)* when, for all objects U of $\mathcal{C}$ and y of YU, there exists a cover $e : V \longrightarrow U$, an object x of XV, and an isomorphism $(Ye)y \cong f_V x$.

A morphism $f : X \longrightarrow Y$ in $\mathcal{F}$ is called *cover cartesian fully faithful (c.c.f.f.)* when it is (pointwise) f.f. and, for all covers $e : V \longrightarrow U$, the following square is equivalent to a pseudopullback.

$$\begin{array}{ccc} XU & \xrightarrow{f_U} & YU \\ {\scriptstyle Xe}\downarrow & \overset{\cong}{\underset{f_e}{}} & \downarrow{\scriptstyle Ye} \\ XV & \xrightarrow[f_V]{} & YV \end{array}$$

**Proposition 5.2** *The classes of l.s.o. and c.c.f.f. morphisms form a regular factorization system on $\mathcal{F}$ as a bicategory.*

**Proof** Given any morphism $f : X \longrightarrow Y$ in $\mathcal{F}$, define ZU to be the full subcategory of YU consisting of those objects y for which there exists a cover $e : V \longrightarrow U$, an object x of XV, and an isomorphism $(Ye)y \cong f_V x$. Let $i_U$ be the inclusion $i_U : ZU \longrightarrow YU$. Because covers form a calculus of left fractions, we see that, for all $u : W \longrightarrow U$, the functor Yu restricts to a functor $Zu : ZU \longrightarrow ZW$; so Z is an object of $\mathcal{F}$ and $i : Z \longrightarrow Y$ is f.f. (and strict in fact). It is readily checked that $i$ is indeed cover cartesian. Also, since each $f_U$ lands in $ZU$, we obtain the components of a morphism $t : X \longrightarrow Z$ with $f = i\,t$. It is clear that $t$ is l.s.o. The remaining details are routine. **Q.E.D.**

## 6. Classification of locally trivial structures

We require four ingredients:

(a) a category $\mathcal{C}$;

(b) a category $X$ parametrized by $\mathcal{C}$;

(c) a *cover* $e : V \longrightarrow U$;

(d) a family $t$ of *trivial objects* of $X$.

More explicitly, $\mathcal{C}$ can be any finitely complete category, $X$ can be any pseudofunctor $X : \mathcal{C}^{op} \longrightarrow \mathrm{Cat}$, $e : V \longrightarrow U$ can be any morphism of $\mathcal{C}$, and $t$ is an object of $XT$ for some object $T$ of $\mathcal{C}$. We think of $t$ as a family of objects of type $X$ parametrized by $T$; sometimes we identify it with the corresponding morphism $t : T \longrightarrow X$ in $\mathcal{F}$. Localizing will be understood with respect to the view that our morphism $e : V \longrightarrow U$ is a cover.

Let $\mathrm{Loc}(t\,;e)$ be the full subcategory of $XU$ consisting of the objects $x$ for which there exist a morphism $z : V \longrightarrow T$ in $\mathcal{C}$ and an isomorphism $(Xe)x \cong (Xz)t$. So the objects of $\mathrm{Loc}(t\,;e)$ are thought of as U-families of objects of type $X$ that are locally isomorphic to trivial objects. A more bicategorical definition of $\mathrm{Loc}(t\,;e)$ is as follows. Let $Q(t\,;e)$ denote the category obtained as the following pseudopullback.

$$\begin{array}{ccc} Q(t\,;e) & \xrightarrow{q} & \mathcal{F}(V,T) \\ {\scriptstyle p}\downarrow & \cong & \downarrow{\scriptstyle t_V} \\ XU & \xrightarrow[Xe]{} & XV \end{array}$$

Then $p$ factors as a composite

$$Q(t\,;e) \xrightarrow{p_1} \mathrm{Loc}(t\,;e) \xrightarrow{j_1} XU$$

where $p_1$ is e.s.o. and $j_1$ is the f.f. inclusion.

A factorization

$$\begin{array}{ccccc} V & & \xrightarrow{e} & & U \\ & {\scriptstyle s}\searrow & \cong & \nearrow{\scriptstyle e'} & \\ & & P & & \end{array}$$

in $\mathcal{F}$ is said to be *of effective descent for* X when s is e.s.o. and $e'$ induces an equivalence of categories $XU \simeq \mathcal{F}(P,X)$. It is expected that P should be in Cat$\mathcal{C}$, but that is not really necessary.

Define X[t] by factoring in $\mathcal{F}$ as follows:

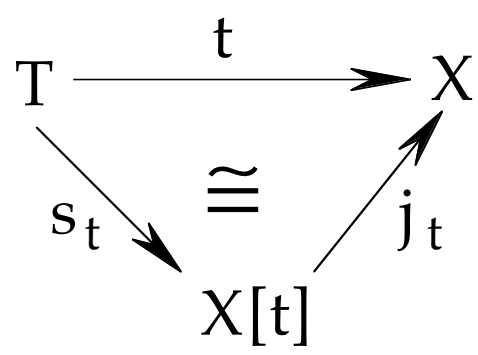

where $s_t$ is e.s.o. and $j_t$ is f.f. Again, it is expected that X[t] should be in Cat$\mathcal{C}$ (and there is some chance of this when X has small homs), but again this is not really necessary.

The following result is essentially from [JSS] and contains the categorical version Fundamental Theorem of Galois Theory due to [Jdz].

**Theorem 6** *There is an equivalence of categories*

$$\mathcal{F}(P,X[t]) \simeq \mathrm{Loc}(t\,;e)\,.$$

**Proof** By the bicategorical factorization system property, since $s : V \longrightarrow P$ is e.s.o. and $j_t : X[t] \longrightarrow X$ is f.f., the bottom right square below is equivalent to a pseudopullback.

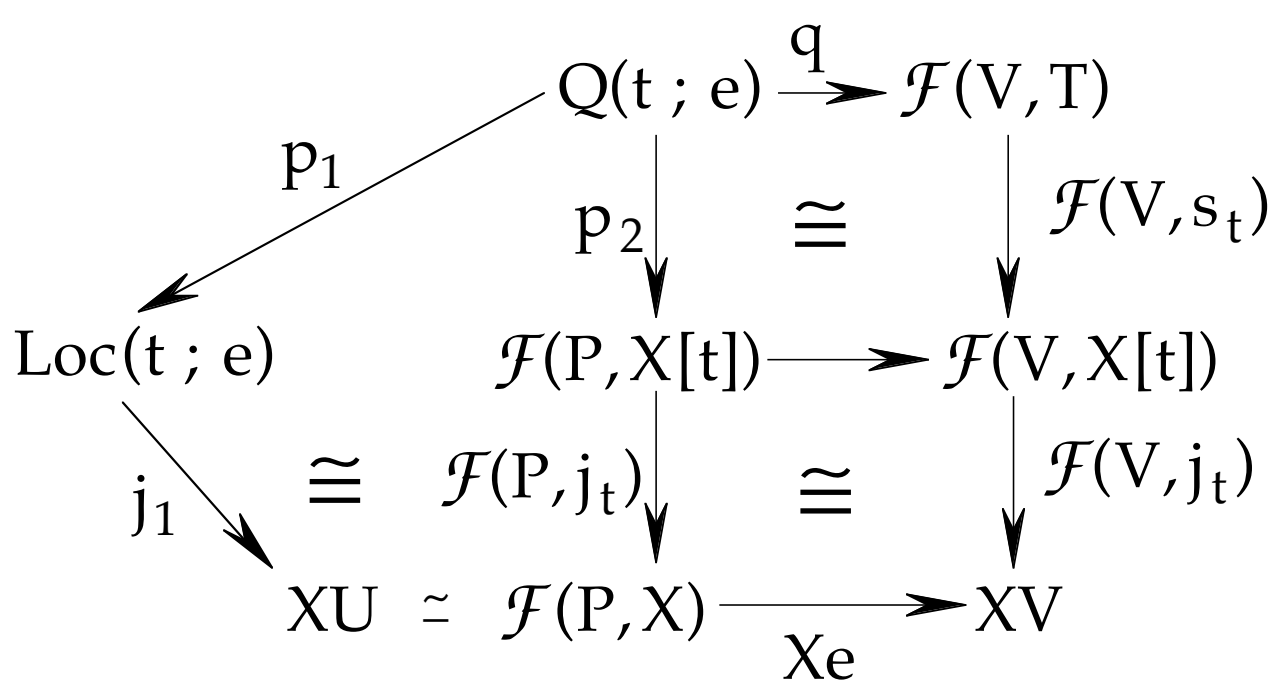

It follows from the definition of Q(t ; e) that there exists a functor $p_2$ as in the above diagram such that the top right square is equivalent to a pseudopullback. Since V is in $\mathcal{C}$ and $s_t$ is e.s.o., the functor $\mathcal{F}(V,s_t)$ is e.s.o. It follows by regularity that $p_2$ is e.s.o. Thus the left-hand region of the diagram provides two factorizations of p into an e.s.o. and an f.f. The images are therefore equivalent. **Q.E.D.**

By way of a typical example, take $\mathcal{C}$ to be the category Top of topological spaces. Define the pseudofunctor $X : \mathrm{Top}^{\mathrm{op}} \longrightarrow \mathrm{Cat}$ to take a space U to the category XU of modules in Top / U over the ring object $\mathrm{pr}_2 : \mathbf{R} \times U \longrightarrow U$ where $\mathbf{R}$ is the topological ring of real numbers. Let K be the space of pairs $(n, x)$ where n is a natural number and x is a vector in n-dimensional real Euclidean space. Then the first projection $K \longrightarrow \mathbf{N}$ is a module over $\mathbf{R} \times \mathbf{N} \longrightarrow \mathbf{N}$, where $\mathbf{N}$ is the discrete space of natural numbers, and so is an object t of $X\mathbf{N}$. Let $\mathcal{U} = (U_i)$ be an open cover of the space U and let V be the disjoint union $V = \sum_i U_i$ with $e : V \longrightarrow U$ induced by the inclusions $U_i \longrightarrow U$. Then X[t] can be taken to be the topological category Mat(**R**) whose objects are natural numbers and whose morphisms $n \longrightarrow m$ are $m \times n$ matrices. The topological category called P above is none other than the nerve $\mathrm{Ner}\,\mathcal{U}$ of the covering $\mathcal{U}$. Theorem 6 gives an equivalence between the category

$$\mathrm{Cat(Top)}(\mathrm{Ner}\,\mathcal{U}, \mathrm{Mat}(\mathbf{R}))$$

of topological functors from the nerve of $\mathcal{U}$ to Mat(**R**) and the category of real vector bundles over U trivialized by the covering $\mathcal{U}$. This yields the clutching constructions for vector bundles and, on restricting the equivalence to the groupoids of invertible morphisms, yields the classification of vector bundles by Čech 1-cocycles with coefficients in the real general linear groups $GL_n(\mathbf{R})$.

## 7. Stacks and torsors

Suppose (as near the end of Section 5) we have a finitely complete category $\mathcal{C}$ with a calculus of left fractions whose morphisms are called *covers*. For each cover $e : V \longrightarrow U$, we can form the category Er(e) in $\mathcal{C}$ called the *equivalence relation* for e : it is the simplicial object

$$\cdots\ V \times_U V \times_U V \;\substack{\longrightarrow \\ \longrightarrow \\ \longrightarrow}\; V \times_U V \;\substack{\xrightarrow{pr_2} \\ \xleftarrow{diag} \\ \xrightarrow{pr_1}}\; V\ .$$

We have a factorization

$$\begin{array}{ccc} V & \xrightarrow{\quad e \quad} & U \\ & {}_{s}\searrow \quad \nearrow{}_{e'} & \\ & \mathrm{Er}(e) & \end{array}$$

in which s is b.o. An object X of $\mathcal{F}$ is said to be *1-separated* when the functor

$$\mathcal{F}(e', X) : \mathcal{F}(U, X) \longrightarrow \mathcal{F}(\mathrm{Er}(e), X)$$

is faithful for all covers e. The object X is said to be *2-separated* when the displayed functor is fully faithful for all covers e. We call X a *stack* when $\mathcal{F}(e',X)$ is an equivalence of categories for all covers e.

In other words, an object X of $\mathcal{F}$ is a stack (for the given covers in $\mathcal{C}$ ) when, for all covers e, the above displayed factorization of e is of effective descent for X.

For any $t: T \longrightarrow X$ with T in $\mathcal{C}$, we put

$$\mathrm{Loc}_X(t)U = \bigcup_{e \in \mathrm{Cov}U} \mathrm{Loc}(t;e).$$

Then, $\mathrm{Loc}_X(t)$ becomes an object of $\mathcal{F}$; indeed, it is the l.s.o./c.c.f.f. image of $t: T \longrightarrow X$. Theorem 6 yields the equivalence

$$\mathrm{Loc}_X(t)U \simeq \operatorname*{colim}_{e \in \mathrm{Cov}U} (\mathrm{Cat}\mathcal{C})(\mathrm{Er}(e), X[t]),$$

where the right-hand side is a filtered colimit in Cat (and so commutes with finite limits)

Let A be a category in $\mathcal{C}$. An A-*torsor trivialized by a cover* $e: V \longrightarrow U$ is a discrete fibration E from A to U for which there exist a morphism $a: V \longrightarrow A$ and a commutative diagram in which the square is a pullback. In other words, A-torsors trivialized by e are the objects of the category Loc(t ; e) where t is the composite

$$A_0 \longrightarrow A \xrightarrow{y_A} \mathcal{P}A;$$

put $\mathrm{Tors}(A;e) = \mathrm{Loc}(t;e)$ for this t. So "trivial" here means "representable" in the sense of being in the image of the yoneda morphism. An A-*torsor at* U is an A-torsor trivialized by some cover $e: V \longrightarrow U$. We put $\mathrm{Tors}A = \mathrm{Loc}_{\mathcal{P}A}(t)$, an object of $\mathcal{F}$.

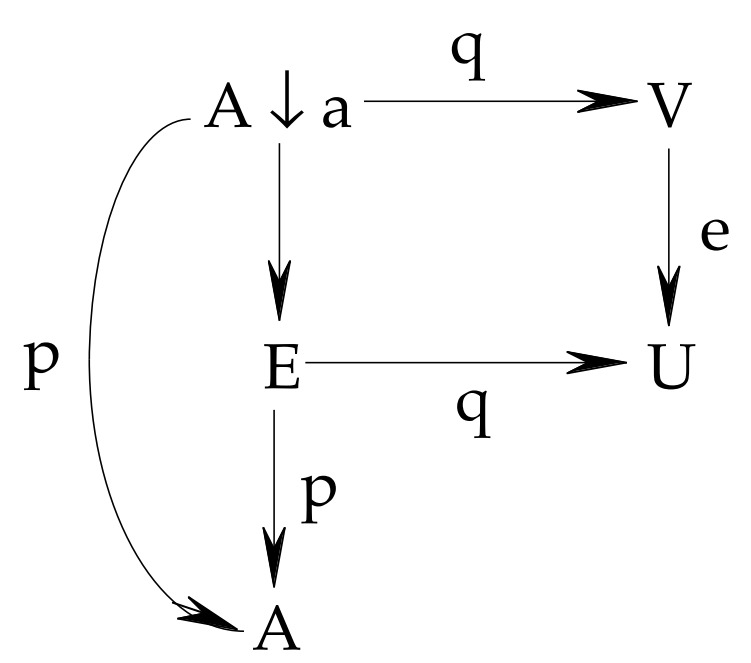

As a corollary of Theorem 6 we have the equivalence of categories

$$\mathrm{Tors}(A;e) \simeq (\mathrm{Cat}\mathcal{C})(\mathrm{Er}(e), A)$$

and the equivalence

$$\mathrm{Tors}A \simeq \operatorname*{colim}_{e \in \mathrm{Cov}U} (\mathrm{Cat}\mathcal{C})(\mathrm{Er}(e), A)$$

in $\mathcal{F}$ which interpret as saying that *all* A-*torsors can be constructed from* $\hat{\text{C}}$*ech cocycles with coefficients in* A.

Now we point out the fundamental relationship between stacks and torsors (see [St3] and [St4]).

**Theorem 7.1** *An object of $\mathcal{F}$ is a stack if and only if it admits all colimits weighted by torsors. These colimits are absolute: that is, preserved by all morphisms in $\mathcal{F}$.*

**Proof** Suppose X is a stack. Take any torsor E from A to U and $f : A \longrightarrow X$. Let $e : V \longrightarrow U$ be a cover trivializing E and let $a : Er(e) \longrightarrow A$ be the "cocycle" corresponding to E. Then we have $fa$ in the category $\mathcal{F}(Er(e), X)$ which is equivalent to XU since X is a stack. The object of XU corresponding to $fa$ is col(E,f). Conversely, suppose X is cocomplete with respect to torsors as weights. We need to prove that the functor $\mathcal{F}(U,X) \longrightarrow \mathcal{F}(Er(e),X)$ induced by $e' : Er(e) \longrightarrow U$ is an equivalence. We use the fact that $e' \downarrow U$ is an Er(e)-torsor trivialized by e. The required inverse equivalence is defined by the colimit $\text{col}(e' \downarrow U, -)$ weighted by $e' \downarrow U$.

For the second sentence of the Theorem, it suffices to show that the colimit is preserved by any morphism $h : X \longrightarrow Y$ into a cocomplete object Y. In particular, Y is a stack. So we see that col(E, h f) in YU corresponds to $hfa$. By evaluating the following commutative square at col(E,f), we obtain the isomorphism $\text{col}(E, hf) \cong h\,\text{col}(E, f)$.

$$\begin{array}{ccc} \mathcal{F}(U,X) & \xrightarrow{\mathcal{F}(e',X)} & \mathcal{F}(Er(e),X) \\ {\scriptstyle \mathcal{F}(U,h)}\downarrow & & \downarrow{\scriptstyle \mathcal{F}(Er(e),h)} \\ \mathcal{F}(U,Y) & \xrightarrow[\mathcal{F}(e',Y)]{} & \mathcal{F}(Er(e),Y) \end{array}$$

**Q.E.D.**

Constructing the *associated stack* of an arbitrary object P of $\mathcal{F}$ is therefore the cocompletion of P with respect to torsors. This can be done in various ways. The approach that is closest to the original associated sheaf construction described by Grothendieck [An] is to define LP in $\mathcal{F}$ by

$$(LP)U = \underset{e \in CovU}{\text{colim}}\, \mathcal{F}\big(Er(e), P\big).$$

The morphisms $e' : Er(e) \longrightarrow U$ induce

$$PU \xrightarrow{\sim} \mathcal{F}(U, P) \xrightarrow{\mathcal{F}(e', P)} \mathcal{F}\big(Er(e), P\big),$$

and thereby a morphism $\eta : P \longrightarrow LP$ in $\mathcal{F}$. The proof of the following result can essentially be found in [St2]. It proceeds in three steps by showing that P in $\mathcal{F}$ is 1-separated iff $\eta$ is (pointwise) faithful, that P is 2-separated iff $\eta$ is fully faithful, and that P is a stack iff $\eta$ is an eqivalence.

**Theorem 7.2** *If* P *is any object of* $\mathcal{F}$ *then* $L^3P$ *is the associated stack of* P *in the sense that* $L^3P$ *is a stack and, for all stacks* X *in* $\mathcal{F}$, *the morphism* $P \longrightarrow L^3P$, *obtained by composing three instances of* $\eta$, *induces an equivalence of categories*

$$\mathcal{F}\left(L^3P, X\right) \simeq \mathcal{F}(P, X).$$

## 8. Parity complexes

Free categories on circuit-free directed graphs have particularly simple descriptions. We generalise this to higher dimensions following [St7].

A parity complex C of dimension n consists of a graded set $C = \sum_{0 \le k \le n} C_k$ and functions $(-)^-$ and $(-)^+ : C_k \longrightarrow \mathcal{P}C_{k-1}$ for $0<k\le n$, where $\mathcal{P}$ denotes the power set. For any subset S of $C_k$, we write $S^-$ for the subset of $C_{k-1}$ consisting of all elements in some $x^-$ with $x \in S$; similarly define $S^+$. There are some axioms such as

$$x^- \cap x^+ = \varnothing \qquad \text{and} \qquad x^{--} \cup x^{++} = x^{-+} \cup x^{+-}.$$

The *solid triangle order* $\blacktriangleleft$ on the set C is defined to be the smallest reflexive transitive relation having $x \blacktriangleleft y$ when either $x \in y^-$ or $y \in x^+$. A strong axiom of loop freeness on a parity complex is that *the solid triangle order should be antisymmetric* ; moreover, for the important examples of simplexes, cubes and globes defined below, the order is *linear* (that is, total).

The model for the free n-category $\mathcal{O}C$ on C will now be succinctly described in a purely combinatorial way. An n-cell of $\mathcal{O}C$ is a pair $(M, P)$ of non-empty finite subsets M (for "minus") and P (for "plus") of C such that the following conditions hold (where $\neg S$ means the complement of S in C):

(i) each of M and P contains at most one element of $C_0$ and, for all $x \neq y$ in $C_k$ with $k>0$, if both $x, y \in M$ or if both $x, y \in P$, then the set $\left(x^- \cap y^-\right) \cup \left(x^+ \cap y^+\right)$ is empty;

(ii) $$P = \left(M \cup M^+\right) \cap \neg M^-, \qquad M = \left(P \cup M^-\right) \cap \neg M^+,$$
$$P = \left(M \cup P^+\right) \cap \neg P^-, \qquad M = \left(P \cup P^-\right) \cap \neg P^+.$$

The *k-source* and *k-target* of $(M,P)$ are defined as follows (where $S_k = C_k \cap S$ and $S^{(k)} = \sum_{h \le k} S_h$ for any subset $S$ of $C$):

$$s_k(M,P) = \left(M^{(k)}, M_k \cup P^{(k-1)}\right), \; t_k(M,P) = \left(M^{(k-1)} \cup P_k, P^{(k)}\right).$$

An ordered pair of cells $(M,P)$, $(N,Q)$ is called *k-composable* when

$$t_k(M,P) = s_k(N,Q),$$

in which case their *k-composite* is defined by

$$(M,P) \circ_k (N,Q) = \left(M \cup (N \cap \neg N_k), (P \cap \neg P_k) \cup Q\right).$$

The *k-cells* of $\mathcal{O}C$ are the n-cells $(M,P)$ with $s_k(M,P) = (M,P)$. The proof that $\mathcal{O}C$ is an n-category is non-trivial (and requires more axioms on the parity complex than those mentioned above). There is a dimension preserving injective function

$$x \mapsto \langle x \rangle \; : C \longrightarrow \mathcal{O}C$$

given inductively as follows: for $x \in C_k$, put $\langle x \rangle = (M,P)$ where

$$M_k = P_k = \{x\},$$

$$M_{r-1} = (M_r)^- \cap \neg (M_r)^+, \quad \text{and} \quad P_{r-1} = (P_r)^+ \cap \neg (P_r)^- \quad \text{for } 0 < r \le k.$$

My notation for this particular $M$ and $P$ is $\mu(x)$ and $\pi(x)$ so that $\langle x \rangle = (\mu(x), \pi(x))$. It is also non-trivial to prove that $\mathcal{O}C$ is the *free* n-category generated by the cells $\langle x \rangle$, $x \in C$.

The *product* $C \times D$ of two parity complexes $C, D$ is given by

$$(C \times D)_n = \sum_{p+q=n} C_p \times D_q, \quad (x,a)^\varepsilon = x^\varepsilon \times \{a\} \cup \{x\} \times a^{\varepsilon(p)}$$

for $x \in C_p$, $a \in D_q$ and $\varepsilon \in \{-,+\}$ where $\varepsilon(p) \in \{-,+\}$ is $\varepsilon$ for $p$ even and is not $\varepsilon$ for $p$ odd.

Parity complexes can be regarded as combinatorial chain complexes. Each parity complex $C$ gives rise to a chain complex $FC$ by taking the free abelian groups on each $C_n$ and using the differential $d(x) = x^+ - x^-$, where we have identified $x^+$ with the formal sum of its elements. It is easy to see that we have a canonical isomorphism of chain complexes:

$$F(C \times D) \cong FC \otimes FD,$$

where we remind readers that the tensor-product boundary formula is

$$d(x \otimes a) = dx \otimes a + (-1)^p x \otimes da \qquad \text{for } x \in FC_p \text{ and } a \in FD_q.$$

There are explicit formulas for $\mu(x,a)$ and $\pi(x,a)$ in terms of $\mu(x)$, $\mu(a)$, $\pi(x)$ and $\pi(a)$. To express these, write $\chi^r$ to denote $\chi \in \{\mu,\pi\}$ when $r$ is even and to denote the other element of $\{\mu,\pi\}$ when $r$ is odd. Then

$$\chi(x,a)_n = \bigcup_{r+s=n} \chi(x)_r \times \chi^r(a)_s \,.$$

The *join* $C \bullet D$ of two parity complexes $C$ and $D$ is given by

$$(C \bullet D)_n = C_n + \sum_{p+q+1=n} C_p \times D_q + D_n$$

in which the summands $C$ and $D$ are embedded as parity subcomplexes and the elements $(x,a) \in C_p \times D_q$ are written as $xa$ with

$$(xa)^- = x^- a \cup xa^- \quad \text{and} \quad (xa)^+ = x^+ a \cup xa^+ \quad \text{for } p \text{ odd,}$$

$$(xa)^- = x^- a \cup xa^+ \quad \text{and} \quad (xa)^+ = x^+ a \cup xa^- \quad \text{for } p \text{ even,}$$

where, for example, $x^+a = \{ya : y \in x^+\}$ is taken to mean $\{a\}$ when $p = 0$. In particular, when $D$ consists of a single element $\infty$ in dimension 0, the join $C \bullet D$ is called the *right cone of* $C$ and denoted by $C^>$. Also $D \bullet C$ is the *left cone of* $C$ and denoted by $C^<$.

Let $\mathbf{I}^0$ denote the *parity point;* it is the parity complex $C$ with $C_0 = \{0\}$ and $C_n = \varnothing$ for $n > 0$. The *parity interval* is the parity complex which is the join $\mathbf{I} = \mathbf{I}^0 \bullet \mathbf{I}^0$.

The *parity n-simplex* is the (n+1)-fold join $\mathbf{\Delta}^n = \underbrace{\mathbf{I}^0 \bullet \mathbf{I}^0 \bullet \ldots \bullet \mathbf{I}^0}_{n+1}$ of parity points. In fact, the elements of $\left(\mathbf{\Delta}^n\right)_k$ can be taken to be k-element subsets of $\{0,1,\ldots,n\}$ where $x^-$ consists of the "odd faces" and $x^+$ the "even faces" for such a subset $x$. For $n = 3$:

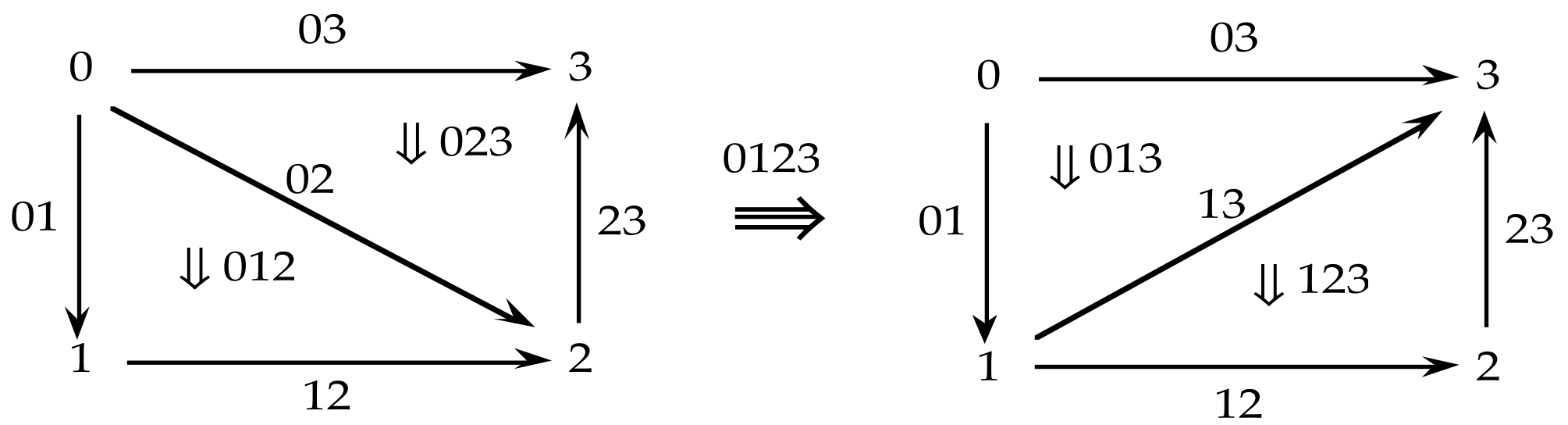

The *parity n-cube* is the n-fold product $\mathbf{I}^n = \underbrace{\mathbf{I} \times \mathbf{I} \times \ldots \times \mathbf{I}}_{n}$ of parity intervals. For $n = 3$:

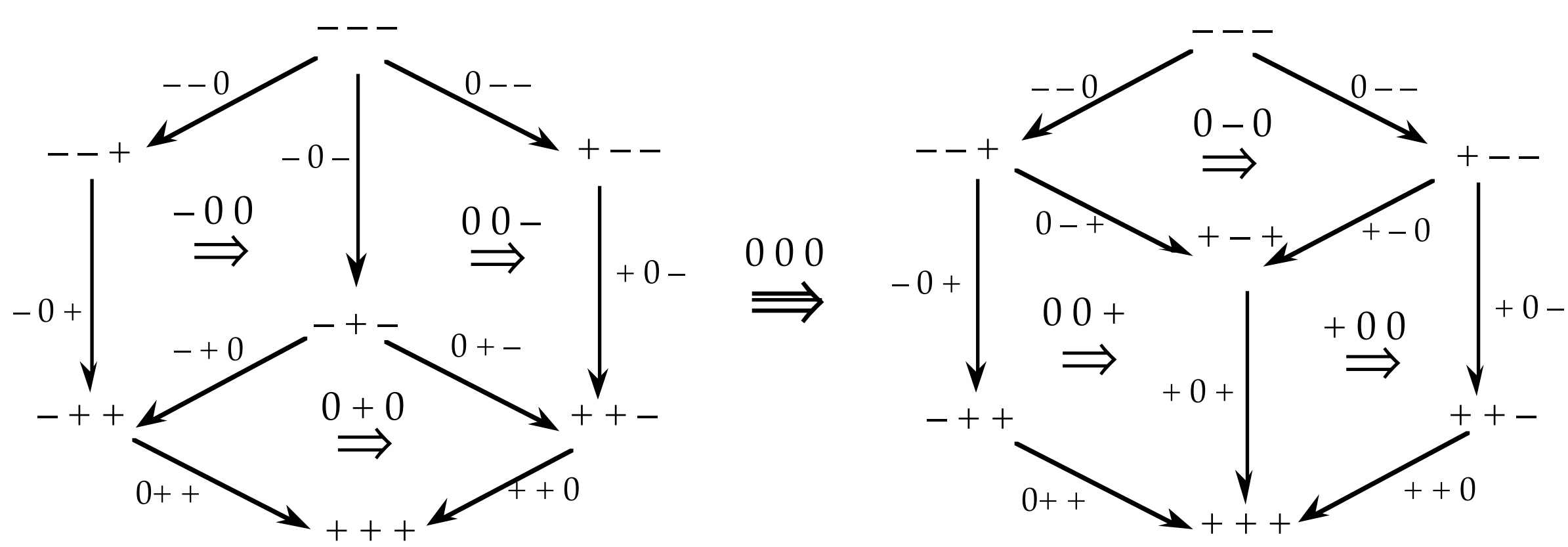

The *parity n-glob* is the parity complex $\mathbf{G}^n$ defined by

$$\left(\mathbf{G}^n\right)_m = \{(\varepsilon, m) : \ \varepsilon = - \text{ or } +\} \quad \text{for } m < n, \qquad \left(\mathbf{G}^n\right)_n = \{n\},$$

$$(\varepsilon,\ m)^- = \{(-, m-1)\}, \quad (\varepsilon,\ m)^+ = \{(+, m-1)\}, \quad n^- = (-, n-1), \quad n^+ = (+, n-1).$$

For $n = 3$:

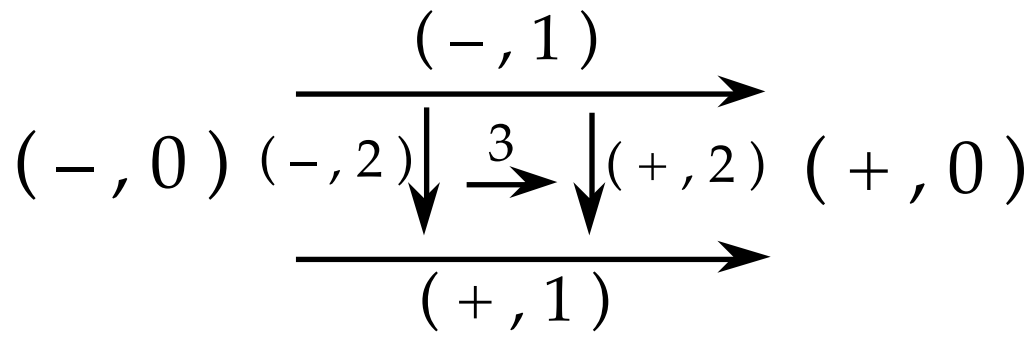

A precise definition of the the free n-category on the n-simplex, called the *n-th oriental,* is

$$\mathcal{O}_n = \mathcal{O}\mathbf{\Delta}^n.$$

A precise definition of the nerve NerA of an ω-category A is then

$$(\mathrm{NerA})_n = \omega\text{-}\mathrm{Cat}(\mathcal{O}_n, A) .$$

This process is quite like Kan's definition of the "singular functor" going from spaces to simplicial sets, so there is also the analogue of a "geometric realization". From the functor $\mathcal{O}_\bullet : \mathbf{\Delta} \longrightarrow \omega\text{-}\mathrm{Cat}$, we obtain the *nerve functor* $\mathrm{Ner} : \omega\text{-}\mathrm{Cat} \longrightarrow \left[\mathbf{\Delta}^{\mathrm{op}}, \mathrm{Set}\right]$ with a left adjoint $\Phi$. While the restriction of Ner to 1-categories is fully faithful, it is not true that Ner itself is full: simplicial maps $\mathrm{NerA} \longrightarrow \mathrm{NerB}$ amount to normal lax functors $A \longrightarrow B$.

## 9. The Gray tensor product of ω-categories and the descent ω-category

We begin by reminding the reader of the technique for left Kan extending monoidal structures along dense functors due to Brian Day [D1], [D2] (whose results more generally cover promonoidal enriched categories).

**Proposition 9** *Suppose* $J : \mathcal{C} \longrightarrow \mathcal{X}$ *is a dense functor from a small monoidal category* $\mathcal{C}$ *into a complete and cocomplete category* $\mathcal{X}$. *The formula*

$$X \otimes Y = \int^{C,D} \left(\mathcal{X}(JC, X) \times \mathcal{X}(JD, Y)\right) \bullet J(C \otimes D)$$

*defines a (left and right) closed monoidal structure on* $\mathcal{X}$ *with* J *strong monoidal if and only if there exist functors* H *and* $H' : \mathcal{C}^{\mathrm{op}} \times \mathcal{X} \longrightarrow \mathcal{X}$ *and isomorphisms*

$$\mathcal{X}\left(JB, H(C,X)\right) \cong \mathcal{X}\left(J(B \otimes C), X\right) \cong \mathcal{X}\left(JC, H'(B,X)\right)$$

*natural in objects* B *and* C *of* $\mathcal{C}$ *and* X *of* $\mathcal{X}$.

For example, when J is the Yoneda embedding of $C$, the tensor product on the presheaf category is *convolution*.

The technique of Proposition 9 was used by the author in [St6] to construct the Gray tensor product of 2-categories. This can be modified to obtain a Gray-like tensor product for ω-categories.

The free ω-categories $\mathcal{O}\mathbf{I}^n$ on the parity cubes (n ≥ 0) form a dense full subcategory $\mathcal{Q}$ of the category $\omega\text{-Cat}$; this essentially amounts to the fact that all possible composites of cells can be found occuring in these cube. The subcategory $\mathcal{Q}$ is monoidal via the tensor product defined by

$$\left(\mathcal{O}\mathbf{I}^m\right)\otimes\left(\mathcal{O}\mathbf{I}^n\right) = \mathcal{O}\mathbf{I}^{m+n}.$$

With some work to satisfy the hypotheses of Proposition 9, we obtain a monoidal structure on $\omega\text{-Cat}$. It is *not* the cartesian monoidal structure. We shall call it the *Gray monoidal structure* on $\omega\text{-Cat}$, although it is not really what John Gray defined; his tensor product was on 2-Cat. The present structure was considered by Richard Steiner [Sn] and explored by Sjoerd Crans [C]. Dominic Verity [V] has another elegant approach using cubical sets. To obtain Gray's original tensor product [Gy1] we need to render all 3-cells identities, although his approach to coherence [Gy2] used the braid groups. To see the connection, consider the braid category $\mathbf{B}$ (as defined in [JS2]) which is the disjoint union of all the usual braid groups as 1-object categories. There is a 2-category $\Sigma\mathbf{B}$ with one object, with hom-category $\mathbf{B}$ and with addition of braids as composition. There is an ω-functor $P : \mathcal{O}\mathbf{I}^\infty \longrightarrow \Sigma\mathbf{B}$ which is universal with the property that it equates all objects, inverts all 2-cells, and takes all 3-cells to identities. Actually, in [St6], the author used the "braid monoids with zero" which are finite monoids that came out of his joint work with Samuel Eilenberg.

Dominic Verity has shown that, for a wide class of parity complexes C, D, we have

$$(\mathcal{O}\mathrm{C})\otimes(\mathcal{O}\mathrm{D}) \cong \mathcal{O}(\mathrm{C}\times\mathrm{D}).$$

Simplexes, cubes, globes, and products of them belong to the class. We shall make use of this result.

To give some further feeling for this Gray tensor product, we shall make a connection with the ordinary tensor product of chain complexes. Each chain complex R gives rise to an ω-category $\vartheta\mathrm{R}$ whose 0-cells are 0-cycles $a \in R_0$, whose 1-cells $b : a \longrightarrow a'$ are elements $b \in R_1$ with $d(b) = a' - a$, whose 2-cells $c : b \longrightarrow b'$ are elements $c \in R_2$ with $d(c) = b' - b$, and so on. All compositions are addition. A functor $\vartheta : \mathrm{DG} \longrightarrow \omega\text{-Cat}$ from the category DG of chain complexes and chain maps. In fact, $\vartheta : \mathrm{DG} \longrightarrow \omega\text{-Cat}$ is a *monoidal functor* where DG has the usual tensor product of chain complexes and $\omega$-Cat has the Gray tensor product. By applying $\vartheta$ on homs, we obtain a (2-) functor $\vartheta_* : \mathrm{DG}\text{-}\mathrm{Cat} \longrightarrow \mathcal{V}_2\text{-}\mathrm{Cat}$ where $\mathcal{V}_2$ is $\omega$-Cat with the Gray tensor product. In particular, since DG is closed, it is a DG-category and we can apply $\vartheta_*$ to it. The $\mathcal{V}_2$-category $\vartheta_*(\mathrm{DG})$ has chain complexes as 0-cells and chain maps as 1-cells; the 2-cells are chain homotopies and the higher cells are higher analogues of chain homotopies. In the next section we shall see the importance of $\mathcal{V}_2$-categories in the homotopy theory of topological spaces, not just the homotopy theory of chain complexes (which is ordinary homological algebra).

We can now solve the problem of defining the descent $\omega$-category of a cosimplicial $\omega$-category. We make considerable use of the fact, mentioned before, that n-categories are models of a finite-limit theory. Such models have their structure preserved by left-exact functors and inherited by representing objects. For example, the functor $\mathrm{Cell}_n : \omega\text{-}\mathrm{Cat} \longrightarrow \mathrm{Set}$, which assigns the set of n-cells to each $\omega$-category, is represented by the free n-category $\mathcal{O}\mathbf{G}^n$ on the n-glob: that is,

$$\mathrm{Cell}_n A \;\cong\; \omega\text{-}\mathrm{Cat}(\mathcal{O}\mathbf{G}^n, A).$$

The set of n-cells in an $\omega$-category forms an n-category; so $\mathcal{O}\mathbf{G}^n$ *is a co-n-category in the category* $\omega$-Cat. Now using the fact that co-n-categories are taken to co-n-categories by right-exact functors, we see that $\mathcal{O}\mathbf{G}^n \otimes A$ is a co-n-category in $\omega$-Cat for all $\omega$-categories A. In particular,

$$\mathcal{O}\mathbf{G}^n \otimes \mathcal{O}_m \;=\; \mathcal{O}\mathbf{G}^n \otimes \mathcal{O}\mathbf{\Delta}^m \;=\; \mathcal{O}\left(\mathbf{G}^n \times \mathbf{\Delta}^m\right)$$

is a co-n-category in $\omega$-Cat for all $m \geq 0$.

Allowing m to vary, we obtain a co-n-category $\mathcal{O}\left(\mathbf{G}^n \times \mathbf{\Delta}^\bullet\right)$ in the category $\left[\mathbf{\Delta}, \omega\text{-}\mathrm{Cat}\right]$ of cosimplicial $\omega$-categories. Hence, for any cosimplicial n-category $\mathcal{E}$, we obtain an n-category

$$\mathrm{Desc}\,\mathcal{E} \;=\; \left[\mathbf{\Delta}, \omega\text{-}\mathrm{Cat}\right]\left(\mathcal{O}\left(\mathbf{G}^n \times \mathbf{\Delta}^\bullet\right), \mathcal{E}\right).$$

We thus have our precise definition of the n-category $\mathrm{Desc}\,\mathcal{E}$ (with somewhat more detail than appears in [St7]).

## 10. Weak n-categories, cohomology and homotopy

There are now many plausible definitions of weak n-category; see [Lr]. For any of these we expect a weak 0-category to be a set, a weak 1-category to be a category, a weak 2-category to be a bicategory in the sense of Bénabou [Bu], and a weak 3-category to be a tricategory in the sense of [GPS]. The definition we wish to concentrate on here is that of Batanin as described in [Bn1] and [Bn2]. The starting point is the category of globular sets (or ω-graphs) and the monad on it whose algebras are ω-categories.

We can approach ω-graphs in the same way we approached ω-categories in the Introduction. For any symmetric monoidal category $\mathcal{V}$, there is a symmetric monoidal category $\mathcal{V}\text{-Gph}$ whose objects are $\mathcal{V}$-graphs; a $\mathcal{V}$-graph G has a set $G_0$ of *vertices* together with, for each ordered pair x, y of vertices, an object $G(x,y)$ of $\mathcal{V}$ (the "object of edges"). Starting with the category Set of sets using cartesian product for the monoidal structure, we can iterate the process $\mathcal{V} \mapsto \mathcal{V}\text{-Gph}$ yielding the following sequence of definitions:

$$\text{Set}, \quad \text{Gph} := \text{Set-Gph}, \quad 2\text{-Gph} := \text{Gph-Gph}, \quad 3\text{-Gph} := (2\text{-Gph})\text{-Gph}, \ \ldots$$

all terms having cartesian product as monoidal structure. Each set can be regarded as a discrete graph (the objects of edges are empty) so there are inclusions

$$\text{Set} \subset \text{Gph} \subset 2\text{-Gph} \subset 3\text{-Gph} \subset \ldots .$$

The union of this chain is the category $\omega\text{-Gph}$ of *ω-graphs* [2]. We define the *n-cells* in an ω-graph just as in ω-categories (see the Introduction); each n-cell has a source (n–1)-cell and a target (n–1)-cell. In this way an ω-graph G can be regarded as having the same kind of structure as a parity complex; it is graded by the dimension of the cells, the sets $x^-$ is the singleton consisting of the source of x and $x^+$ is the singleton consisting of the target of x. It follows that we can define the solid triangle order for ω-graphs. The author [St11] has defined an ω-graph to be a *globular cardinal* when the solid triangle order is linear; also see [MZ].

Under reasonable (co)completeness conditions on $\mathcal{V}$, Wolff [W] showed the forgetful functor $\mathcal{V}\text{-Cat} \longrightarrow \mathcal{V}\text{-Gph}$ to have not only a left adjoint but to be monadic; also see [Bi]. It follows that the forgetful functor $U_n : n\text{-Cat} \longrightarrow n\text{-Gph}$ has a left adjoint $F_n$ for all $0 \le n \le \omega$. Indeed, $U_n$ is also monadic; we write $D_n$ for the monad $U_n F_n$ on $n\text{-Gph}$ generated by the adjunction $F_n \dashv U_n$.

The starting point of Batanin's work was his explicit description of $D_n$. For each n-graph G, the m-cells of $D_n G$ are to be thought of as globular pasting diagrams of k-cells for $k \le m$. Batanin was able to code these globular pasting diagrams in terms of plane trees of height m. Each plane tree t of height m gives rise to an m-graph $t^*$. An ω-graph is a globular cardinal if and only if it is isomorphic to $t^*$ for some plane tree t. The cells of $D_n G$ are n-graph morphisms $t^* \longrightarrow G$.

Weak n-categories are expected to have all the composition operations of the strict n-categories, however, these operations are not expected to be strictly associative or strictly functorial over each other as in a strict n-category. Batanin realized that weak n-categories should also be algebras for some monad $K_n$ on $n\text{-}Gph$. Write $Wk\text{-}n\text{-}Cat$ for the category of Eilenberg-Moore algebras for $K_n$; the objects are weak n-categories but the morphisms are very strict, preserving all the structure precisely. Since every strict n-category should be a particular kind of weak one, there should be a monad morphism $K_n \longrightarrow D_n$ inducing the inclusion $n\text{-}Cat \longrightarrow Wk\text{-}n\text{-}Cat$. The genius of Batanin's approach was the idea, inspired by homotopy theory [3], that $K_n$ should be *contractible* in a suitable sense; indeed, $K_n$ should be the initial contractible monad *with a system of compositions.* This "system" ensured that the composition operations and identities available in an n-category were there in the algebras for $K_n$, while contractibility gave the weak associativity and functoriality.

[2] Something a bit bigger than this union is the category of *globular sets*.

Batanin provided a construction for $K_n$ in [Bn1] and [Bn2]; another arose from [Pn] and [Bn4]. Recent work of Batanin seems to be leading to an explicit combinatorial description of $K_n$ using polyhedra constructed from Joyal's morphisms of Batanin's trees as appearing in [Jl] and [BS]. We shall not need much of this detail here: suffice it to say that, like $D_n$, the endofunctor $K_n$ preserves filtered colimits so that its algebras (the weak n-categories) are also models of a finite-limit theory. This means we can take models in any finitely complete category $\mathcal{C}$; that is, we can speak of *weak n-categories internal to* $\mathcal{C}$.

What is more, if we let A be a weak n-category in $\mathcal{C}$ and let R be a simplicial object of $\mathcal{C}$, we obtain a cosimplicial weak n-category $\mathcal{C}(R, A)$. It is important to realize that the coface and codegeneracy morphisms of $\mathcal{C}(R, A)$ are strict; that is, they are morphisms of $Wk\text{-}n\text{-}Cat$. In the following sections we shall indicate how to define the descent weak n-category of such a cosimplicial weak n-category. Then we define *the cohomology* ***weak*** *n-category* $\mathcal{H}(R, A)$ *of* R *with coefficients in* A by

$$\mathcal{H}(R, A) = \mathrm{Desc}\,\mathcal{C}(R, A).$$

As mentioned in the strict case, the cells of $\mathcal{H}(R, A)$ are *cocycles* of R with coefficients in A. To work with cocycles up to *coboundary* is to work with them up to "equivalence". So we shall briefly discuss equivalence in weak n-categories.

It is well known what is meant for two elements in a set to be equal (= 0-equivalent) and what it means for an arrow in a category to be an isomorphism ( = 1-equivalence). It is also well known what it means for a morphism in a bicategory to be an equivalence ( = 2-equivalence). A 1-cell $f : a \longrightarrow b$ in a tricategory is called a *biequivalence* ( = 3-equivalence) when there exists a 1-cell $g : b \longrightarrow a$ such that $fg$ and $gf$ are both equivalent to identity 1-cells.

Notice that for these kinds of equivalences no use is made of associativity or functoriality of composition. In fact it is possible to define *m-equivalence* in algebras for any monad with a system of compositions; in such algebras, there is a composition of m-cells for all $0 < m \leq n$. The definition of *m-equivalence* is recursive: an n-cell is a 1-equivalence when it is invertible; an m-cell $f : a \longrightarrow b$ is an $(n-m+1)$-equivalence when there exists an m-cell $g : b \longrightarrow a$ with $(n-m)$-equivalences $g f \longrightarrow 1_a$ and $f g \longrightarrow 1_b$.

[3] He introduced a higher dimensional notion of *operad* and expressed contractibility in terms of that.

We can define homotopy sets for any weak n-category A. We define $\pi_0(A)$ to be the set of n-equivalence classes of 0-cells of A. Let a be any 0-cell of A and let AutEq(a) denote the full sub-weak-(n–1)-category of $A(a, a)$ whose 0-cells are the n-equivalences $a \longrightarrow a$. We define the *fundamental group* $\pi_1(A, a)$ to be the set $\pi_0(\mathrm{AutEq(a)})$ equipped with the multiplication induced by composition of 1-cells in A. We recursively define homotopy (abelian) groups $\pi_n(A, a)$, $n > 1$, by

$$\pi_{n+1}(A, a) = \pi_n(\mathrm{AutEq(a)}, 1_a).$$

## 11. Computads, descent and simplicial nerves for weak n-categories

Computads were introduced in [St0] to provide presentations of 2-categories that were more efficient than presentations by 2-graphs. Such a computad is a 2-graph whose 0-cells and 1-cells form the underlying category of a free category on a graph; in other words, we are given a graph together with 2-cells between paths in the graph. These computads were later called 2-computads as the author had need for n-computads for all positive integers n; see [Pr], [St9] and [St10].

More recently, Batanin has defined computads, not just for n-categories, but for any algebraic structure on globular sets; see [Bn3] and [Bn5]. For example, computads for bicategories are not the same as computads for 2-categories; the 2-cells in a computad for bicategories have chosen bracketings for their source and target paths. We shall now explain the general definition in terms similar to the case of computads for (strict) n-categories.

Suppose $T_n$ is a monad on $n\text{-}\mathrm{Gph}$ for each natural number n and let $T_n\text{-}\mathrm{Alg}$ denote the category of Eilenberg-Moore algebras. Let $U_n : T_n\text{-}\mathrm{Alg} \longrightarrow n\text{-}\mathrm{Gph}$ be the underlying functor with left adjoint $F_n$. Let $W_{n-1} : n\text{-}\mathrm{Gph} \longrightarrow (n-1)\text{-}\mathrm{Gph}$ be the functor that forgets about n-cells and let $I_n : (n-1)\text{-}\mathrm{Gph} \longrightarrow n\text{-}\mathrm{Gph}$ be the inclusion; indeed $I_n$ is the fully faithful left adjoint of $W_{n-1}$.

A second sequence of monads $\overline{T}_n$ on $n\text{-}\mathrm{Gph}$ can be constructed from the sequence of monads $T_n$. Define $\overline{T}_{n-1}$ to be the right Kan extension of $W_{n-1} T_n$ along $W_{n-1}$; in fact, $\overline{T}_{n-1} = W_{n-1} T_n I_n$. (For the special case where $T_n$ is the monad for strict n-categories, we have $\overline{T}_n = T_n$. However, when $T_n$ is the monad for weak n-categories, $\overline{T}_1$ assigns the graph of bracketed paths in a graph, whereas $T_1$ assigns the usual graph of paths.) Since $W_{n-1}$ is a monad morphism, it induces a functor $\overline{W}_{n-1} : T_n\text{-}\mathrm{Alg} \longrightarrow \overline{T}_{n-1}\text{-}\mathrm{Alg}$ such that $\overline{U}_{n-1}\,\overline{W}_{n-1} = W_{n-1}\,U_n$, where we put bars overtop data pertaining to the $\overline{T}_n$ to distinguish it from the corresponding data for the $T_n$.

For all sequences of monads $T_n$ on $n\text{-}\mathrm{Gph}$, the category $T_n\text{-}\mathrm{Cpd}$ of *n-computads for* $T_n$*-algebras* is defined inductively along with the functor $V_n : T_n\text{-}\mathrm{Alg} \longrightarrow T_n\text{-}\mathrm{Cpd}$ and its left adjoint $L_n \dashv V_n$. For $n = 0$, $T_n\text{-}\mathrm{Cpd}$ is $T_0\text{-}\mathrm{Alg}$ with $V_0$ and $L_0$ the identity functor. For $n > 0$, the category $T_n\text{-}\mathrm{Cpd}$ is defined by the following pullback of categories and functors.

$$\begin{array}{ccc} T_n\text{-}\mathrm{Cpd} & \xrightarrow{\;P_n\;} & \overline{T}_{n-1}\text{-}\mathrm{Cpd} \\ {\scriptstyle Q_n}\big\downarrow & & \big\downarrow{\scriptstyle \overline{U}_{n-1}\overline{L}_{n-1}} \\ n\text{-}\mathrm{Gph} & \xrightarrow[\;W_{n-1}\;]{} & (n-1)\text{-}\mathrm{Gph} \end{array}$$

A functor $V'_n : T_n\text{-}\mathrm{Alg} \longrightarrow n\text{-}\mathrm{Gph}$ is defined by the following limit diagram of functors and natural transformations

$$\begin{array}{ccc} V'_n & \rightrightarrows & I_n\overline{U}_{n-1}\overline{L}_{n-1}\overline{V}_{n-1}\overline{W}_{n-1} \\ \big\downarrow & & \big\downarrow{\scriptstyle I_n\overline{U}_{n-1}(\mathrm{counit})\overline{W}_{n-1}} \\ U_n & \underset{t}{\overset{s}{\rightrightarrows}} & I_n\overline{U}_{n-1}\overline{W}_{n-1} \end{array}$$

where $s$ and $t$ are the natural transformations whose components assign the source and target $(n-1)$-cells to each n-cell. Notice that $W_{n-1}U_n = \overline{U}_{n-1}\overline{W}_{n-1} = W_{n-1}I_n\overline{U}_{n-1}\overline{W}_{n-1}$ and $W_{n-1}s = W_{n-1}t = 1_{W_{n-1}U_n}$; this implies

$$W_{n-1}V'_n = \overline{U}_{n-1}\overline{L}_{n-1}\overline{V}_{n-1}\overline{W}_{n-1}$$

since $W_{n-1}$ preserves limits. Using the pullback property of $T_n\text{-}\mathrm{Cpd}$, there exists a unique functor $V_n : T_n\text{-}\mathrm{Alg} \longrightarrow n\text{-}\mathrm{Cpd}$ such that $P_nV_n = \overline{V}_{n-1}\overline{W}_{n-1}$ and $Q_nV_n = V'_n$. It is proved in [Bn3] that $V_n$ has a left adjoint $L_n$. This completes the inductive definition.

Just as for ordinary operads, the functor $V_n : T_n\text{-}\mathrm{Alg} \longrightarrow T_n\text{-}\mathrm{Cpd}$ is monadic; again see [Bn3].

The author has long held the view that the orientals should be transferable to contexts other than strict n-categories — to weak n-categories, for example. I am grateful to Michael

Batanin for correcting my naive view of how to do this. He points out that each monad morphism $\theta : T_n \longrightarrow D_n$ induces a functor $\theta^* : D_n\text{-Cpd} \longrightarrow T_n\text{-Cpd}$. Some choice is involved in the definition of $\theta^*$ (such as a splitting of $\theta$ as a mere natural transformation) but all choices are essentially equivalent. The full inductive definition of $\theta^*$ must await another paper, however, the idea is clear enough. Take for example the case where $n = 2$ and $T_n$ is the monad whose algebras are bicategories. Given an ordinary computad H, we must create a computad $\theta^*H$ for bicategories. This is done by choosing a bracketing of each source and target path of each 2-cell of H and making that a single 2-cell of $\theta^*H$. This means that each 2-cell of H leads to only one 2-cell of $\theta^*H$; of course, in the free bicategory on $\theta^*H$ there will be 2-cells between the other bracketings of the source and target paths obtained by using the associativity constraints available in the bicategory.

Start with any parity complex C of dimension n. Form the free n-category $\mathcal{O}C$. Take the underlying computad $V_n\mathcal{O}C$ for strict n-categories (that is, it is a $D_n$-computad). Now we apply the functor $\theta^* : D_n\text{-Cpd} \longrightarrow T_n\text{-Cpd}$ to obtain a $T_n$-computad $\theta^* V_n\mathcal{O}C$

Now we apply the functor $L_n : T_n\text{-Cpd} \longrightarrow T_n\text{-Alg}$ to obtain $\mathcal{O}_T C = L_n\theta^* V_n\mathcal{O}C$. We call $\mathcal{O}_T C$ *the free* $T_n$*-algebra on the parity complex* C.

In particular, for the monad $K_n$ for weak n-categories, we have *the free weak n-category* $\mathcal{O}_K C$ *on the parity complex* C.

One application of this is to the descent construction for weak n-categories. For we now have the cosimplicial weak ω-category $\mathcal{O}_K\left(\mathbf{G}^n \times \mathbf{\Delta}^\bullet\right)$; that is, an object of the functor category $[\mathbf{\Delta}, \text{Wk-}\omega\text{-Cat}]$. We believe it will be possible to show that $\mathcal{O}_K\left(\mathbf{G}^n \times \mathbf{\Delta}^\bullet\right)$ is actually a co-weak-n-category in $[\mathbf{\Delta}, \text{Wk-}\omega\text{-Cat}]$. Then, for any cosimplicial weak-n-category $\mathcal{E}$, we would obtain a weak-n-category

$$\text{Desc}\,\mathcal{E} = [\mathbf{\Delta}, \text{Wk-}\omega\text{-Cat}]\left(\mathcal{O}_K\left(\mathbf{G}^n \times \mathbf{\Delta}^\bullet\right), \mathcal{E}\right).$$

A related application is to obtain the simplicial nerve of a weak ω-category. We might call the weak n-category $\mathcal{O}_K\mathbf{\Delta}^n$ *the* $n^{th}$ *weak oriental.* For any weak ω-category A, define the *nerve* NerA of A to be the simplicial set $\text{Wk-}\omega\text{-Cat}\left(\mathcal{O}_K\mathbf{\Delta}^\bullet, A\right)$.

**Conjecture 11.1** *A simplicial set has the form* NerA *for some weak ω-category* A *if and only if it is a weak ω-category in the sense of* [St12].

As a third application, it seems possible to use the descent construction to produce the weak n-category of weak morphisms from one weak n-category to another. Details will appear elsewhere. For the moment we content ourselves with the following remarks on lax functors.

**Proposition 11.2** *The nerve functor* $\mathrm{Ner} : \mathrm{Wk}\text{-}\omega\text{-}\mathrm{Cat} \longrightarrow \left[\boldsymbol{\Delta}^{\mathrm{op}}, \mathrm{Set}\right]$ *commutes with* $\pi_n$ *for all* $n \geq 0$.

Simplicial maps $f : \mathrm{Ner}A \longrightarrow \mathrm{Ner}B$ are *normal lax functors* between the weak ω-categories $A$ and $B$. (The general lax functors are the face morphisms between the simplicial nerves — they are not required to commute with the degeneracies.) Using the familiar process of replacing a map by an inclusion using a mapping cylinder, we see that each such normal lax functor gives rise to a *long exact homotopy sequence.*

$$\pi_n(A,a) \xrightarrow{f_*} \pi_n(B,f(a)) \longrightarrow \pi_n(f,a) \longrightarrow \pi_{n-1}(A,a) \xrightarrow{f_*} \pi_{n-1}(B,f(a))$$

## 12. Brauer groups

Let $\mathcal{M}$ denote a closed braided monoidal category which is finitely cocomplete. We have in mind that $\mathcal{M}$ is the category of modules over a commutative ring R, or the category of finite dimensional comodules for a quantum group. Consider the bicategory $\mathrm{Alm}\mathcal{M}$ whose objects are monoids (also called "algebras") in $\mathcal{M}$, whose morphism $M : A \longrightarrow B$ are left A- right B-bimodules, and whose 2-cells $f : M \Rightarrow M' : A \longrightarrow B$ are module morphisms $f : M \longrightarrow M'$; vertical composition is composition of functions and horizontal composition of modules $M : A \longrightarrow B$, $N : B \longrightarrow C$ is given by tensor product $M \otimes_B N : A \longrightarrow C$ over $B$ (where $M \otimes_B N$ is the coequalizer of the two arrows from $M \otimes B \otimes N$ to $M \otimes N$ given by the actions of $B$ on $M$ and on $N$).

Since $\mathcal{M}$ is braided, the tensor product $A \otimes B$ of algebras is canonically an algebra. This makes $\mathrm{Alm}\mathcal{M}$ into a monoidal bicategory. Let $\Sigma\mathrm{Alm}\mathcal{M}$ denote the 1-object tricategory whose hom bicategory is $\mathrm{Alm}\mathcal{M}$ and whose composition is tensor product of algebras.

In the particular case of the tricategory $\Sigma\mathrm{Alm}\mathcal{M}$, there it is an easy way to find a 3-equivalent Gray category. First replace $\mathcal{M}$ by an equivalent strict monoidal category (see [JS2]). We identify modules $M : A \longrightarrow B$ with left adjoint functors $\left[A^{\mathrm{op}}, \mathcal{M}\right] \longrightarrow \left[B^{\mathrm{op}}, \mathcal{M}\right]$ where $\left[A^{\mathrm{op}}, \mathcal{M}\right]$ is the category of right A-modules in $\mathcal{M}$. The point is that tensor product $M \otimes_B N$ of modules then becomes composition of functors.

Let $\mathcal{B}r(\mathcal{M})$ denote the sub-Gray-category of $\Sigma\mathrm{Alm}\mathcal{M}$ consisting of the arrows A which are biequivalences, the 2-cells M which are equivalences, and the 3-cells f which are isomorphisms. The morphisms A of $\mathcal{Br}(\mathcal{M})$ are called *Azumaya algebras* in $\mathcal{M}$. The 2-cells M of $\mathcal{Br}(\mathcal{M})$ are called *Morita equivalences* in $\mathcal{M}$.

We can form the nerve $\mathrm{Ner}\mathcal{Br}(\mathcal{M})$ of $\mathcal{Br}(\mathcal{M})$. It is a simplicial set whose homotopy objects are of special importance. In particular, $\pi_0\mathrm{Ner}\mathcal{Br}(\mathcal{M})$ is a singleton set, $\pi_1\mathrm{Ner}\mathcal{Br}(\mathcal{M})$ is called the *Brauer group* $\mathrm{Br}(\mathcal{M})$ *of* $\mathcal{M}$, and $\pi_2\mathrm{Ner}\mathcal{Br}(\mathcal{M})$ is the *Picard group* $\mathrm{Pic}(\mathcal{M})$ *of* $\mathcal{M}$. If $\mathcal{M}$ is equivalent to $\mathrm{Mod}(R)$ for a commutative ring R, these are the usual Brauer and Picard groups of R; also $\pi_3\mathrm{Ner}\mathcal{Br}(\mathcal{M})$ is then isomorphic to the group $\upsilon(R)$ of units of R. Compare the approach of Duskin [Dn1].

Now suppose $F : \mathcal{M} \longrightarrow \mathcal{N}$ is a right-exact braided strong-monoidal functor between finitely cocomplete closed braided monoidal categories. (We have in mind the functor $\mathrm{Mod}(\phi) : \mathrm{Mod}(R) \longrightarrow \mathrm{Mod}(S)$ induced by a commutative ring homomorphism $\phi : R \longrightarrow S$.) Such an F determines a weak morphism (compositions are preserved up to equivalence) of tricategories $\mathrm{Alm}F : \mathrm{Alm}\mathcal{M} \longrightarrow \mathrm{Alm}\mathcal{N}$. Weak morphisms preserve n-equivalence for all n. So a weak morphism $\mathcal{Br}(F) : \mathcal{Br}(\mathcal{M}) \longrightarrow \mathcal{Br}(\mathcal{N})$ is induced, and a simplicial map $\mathrm{Ner}\mathcal{Br}(F) : \mathrm{Ner}\mathcal{Br}(\mathcal{M}) \longrightarrow \mathrm{Ner}\mathcal{Br}(\mathcal{N})$ is induced. This leads to the nine term exact sequence

$$\begin{aligned} 1 \longrightarrow \mathrm{Aut}(I_{\mathcal{M}}) \xrightarrow{F_*} \mathrm{Aut}(I_{\mathcal{N}}) \longrightarrow \mathrm{Aut}(F) \longrightarrow \mathrm{Pic}(\mathcal{M}) \xrightarrow{F_*} \mathrm{Pic}(\mathcal{N}) \\ \longrightarrow \mathrm{Pic}(F) \longrightarrow \mathrm{Br}(\mathcal{M}) \xrightarrow{F_*} \mathrm{Br}(\mathcal{N}) \longrightarrow \mathrm{Br}(F) \longrightarrow 1 \end{aligned}$$

in which $\mathrm{Aut}\left(I_{\mathcal{M}}\right)$ denotes the abelian group of automorphisms of the unit $I_{\mathcal{M}}$ for the tensor product in $\mathcal{M}$. Compare with [DI] when $\mathcal{M} = \mathrm{Mod}(R)$.

## §13. Giraud's $H^2$ and the pursuit of stacks

We use Duskin's [Dn2] amelioration of Giraud's theory [Gd2] to show that Giraud's $H^2$ really fits into our general setting for cohomology. We work in a topos $\mathcal{E}$.

A groupoid B in $\mathcal{E}$ is *connected* when $\pi_0 B \cong 1$.

**Lemma 13.1** *Locally connected implies connected.*

**Proof** If $R \longrightarrow 1$ is an epimorphism ("a cover") then the functor $R \times - : \mathcal{E} \longrightarrow \mathcal{E}/R$ reflects isomorphisms (that is, is conservative), and preserves terminal objects and coequalizers. Hence it also reflects coequalizers. So, to see whether

$$B_1 \rightrightarrows B_0 \longrightarrow 1$$

is a coequalizer in $\mathcal{E}$, it suffices to see that

$$R \times B_1 \rightrightarrows R \times B_0 \longrightarrow R$$

is a coequalizer in $\mathcal{E}/R$. **Q.E.D.**

A functor $f : A \longrightarrow B$ in $\mathcal{E}$ is called *e.s.o.* (essentially surjective on objects, as before) when the top composite of $q$ and $d_1$ in the diagram below is an epimorphism $P \longrightarrow B_0$ and the square is a pullback (here **I** is the category with two objects and an isomorphism between them).

$$\begin{array}{ccccc} P & \xrightarrow{\ q\ } & B_0^{\mathbf{I}} & \xrightarrow{\ d_1\ } & B_0 \\ {\scriptstyle p}\downarrow & & \downarrow{\scriptstyle d_0} & & \\ A_0 & \xrightarrow[\ f_0\ ]{} & B_0 & & \end{array}$$

A groupoid $B$ is called a *weak group* when there exists an e.s.o. $b : 1 \longrightarrow B$. In this case, if $G$ denotes the full image of $b$, we have a weak equivalence (that is, e.s.o. fully faithful functor) $G \longrightarrow B$ where $G$ is a group.

**Lemma 13.2** *A groupoid is connected iff it is a locally weak group.*

**Proof** By Lemma 13.2, "if" will follow from "weak group implies connected". Suppose $b : 1 \longrightarrow B$ is e.s.o.; form the pullback $P$ as above with $A = 1$ and $f = b$. To prove

$$B_1 \underset{d_1}{\overset{d_0}{\rightrightarrows}} B_0 \xrightarrow{\ t\ } 1$$

is a coequalizer, suppose $h : B_0 \longrightarrow X$ has $h\,d_0 = h\,d_1$. Then

$$h\,d_1\,q = h\,d_0\,q = h\,b\,p = h\,b\,t\,d_1\,q$$

implies $h = h\,b\,t$ since $d_1\,q$ is epimorphic. So $h$ factors through $t$. However $t$ is a retraction (split by $b$), so the factorization is unique.

Conversely, assume $B$ is connected. Certainly $B_0 \longrightarrow X$ is epimorphic, so we pass to $\mathcal{E}/B_0$ where we pick up a global object $\Delta : B_0 \longrightarrow B_0 \times B$ over $B_0$ which we will see is e.s.o.

$$
\begin{array}{ccccc}
B_1 & \longrightarrow & B_0 \times B_1 & \xrightarrow{1\times d_1} & B_0 \\
{\scriptstyle d_0}\downarrow & & \downarrow{\scriptstyle 1\times d_0} & & \\
B_0 & \xrightarrow[\Delta]{} & B_0 \times B_0 & &
\end{array}
$$

What we must see then is that $(d_0, d_1) : B_1 \longrightarrow B_0 \times B_0$ is epimorphic. Take the epi./mono. factorization of $(d_0, d_1)$ and let K be the image. Since B is a groupoid, K is an equivalence relation on $B_0$. Since $\mathcal{E}$ is exact, K is a kernel pair of its coequalizer. The coequalizer is 1 since B is connected. So the kernel pair is $B_0 \times B_0$. **Q.E.D.**

Recall that the category of groups in a category with finite products is actually a 2-category since group homomorphisms can be regarded as functors; so there are 2-cells amounting to natural transformations. (In fact, we can make it a 3-category by taking central elements of the target group as 3-cells, but this will not be needed here.) So we have a 2-functor

$$\mathrm{Gp} : \mathrm{Cat}_{\times} \longrightarrow 2\text{-}\mathrm{Cat}$$

from the 2-category $\mathrm{Cat}_{\times}$ of categories with finite products and product-preserving functors.

There is a pseudofunctor $\mathcal{E}/- : \mathcal{E}^{\mathrm{op}} \longrightarrow \mathrm{Cat}$ taking an object X of $\mathcal{E}$ to the slice category $\mathcal{E}/X$ and given on morphisms by pulling back along the morphism. It is easy to find an actual 2-functor $\mathbf{E} : \mathcal{E}^{\mathrm{op}} \longrightarrow \mathrm{Cat}$ equivalent to $\mathcal{E}/-$. The composite 2-functor

$$\mathcal{E}^{\mathrm{op}} \xrightarrow{\mathbf{E}} \mathrm{Cat}_{\times} \xrightarrow{\mathrm{Gp}} 2-\mathrm{Cat}$$

defines a 2-category $\mathcal{G}$ in the presheaf category $[\mathcal{E}^{\mathrm{op}}, \mathrm{Set}]$.

It is natural then to look at the cohomology 2-category $\mathcal{H}(\mathcal{E}, \mathcal{G})$ of $\mathcal{E}$ with coefficients in $\mathcal{G}$. What I mean by this is the colimit of all the 2-categories $\mathcal{H}(R, \mathcal{G})$ over all hypercovers R in $\mathcal{E}$, which we regard, via the Yoneda embedding, as special simplicial objects in the category $[\mathcal{E}^{\mathrm{op}}, \mathrm{Set}]$.

What Giraud actually looks at is obtained from $\mathcal{H}(\mathcal{E}, \mathcal{G})$ by lots of quotienting. First form the composite 2-functor

$$\mathcal{E}^{\mathrm{op}} \xrightarrow{\mathcal{G}} 2-\mathrm{Cat} \xrightarrow{\pi_{0*}} \mathrm{Cat}$$

where $\pi_{0*}$ is the 2-functor which applies $\pi_0$ to the hom categories of each 2-category. Let $\mathcal{L} : \mathcal{E}^{\mathrm{op}} \longrightarrow \mathrm{Cat}$ denote the associated stack of that composite 2-functor. The category $\mathcal{L}(X)$ is called *the category of* X-*liens of* $\mathcal{E}$; in particular, $\mathcal{L}(1)$ is the category of *liens* of $\mathcal{E}$.

The stack condition implies that each epimorphism $R \longrightarrow 1$ induces an equivalence between the category $\mathcal{L}(1)$ of liens and the descent category of the following truncated cosimplicial category.

$$\mathcal{L}(R) \; \substack{\longrightarrow \\ \longleftarrow \\ \longrightarrow} \; \mathcal{L}(R \times R) \; \substack{\longrightarrow \\ \longrightarrow \\ \longrightarrow} \; \mathcal{L}(R \times R \times R) \; .$$

Each connected groupoid B determines a lien $\text{lienB} \in \mathcal{L}(1)$ as follows. By Lemma 13.2, there exists an epimorphism $R \longrightarrow 1$ and $G \in \pi_{0*}\mathcal{G}(R)$. The quotient functor $\pi_{0*}\mathcal{G}(R) \longrightarrow \mathcal{L}(R)$ gives an R-lien $[G] \in \mathcal{L}(R)$ which can be enriched with descent data. These descent data are determined up to isomorphism by B. It follows that there is a lien $\text{lienB} \in \mathcal{L}(1)$ taken to B by the functor $\mathcal{L}(1) \longrightarrow \mathcal{L}(R)$.

For any lien L, let $\mathcal{H}^2(\mathcal{E}, L)$ denote the category whose objects are connected groupoids B with $\text{lienB} \cong L$, and whose arrows are weak equivalences of groupoids. We leave as a future quest the study of the 2-category $\mathcal{H}(\mathcal{E}, \mathcal{G})$ versus the categories $\mathcal{H}^2(\mathcal{E}, L)$.

Centre of Australian Category Theory
Mathematics Department, Macquarie Universit